\documentclass[11 pt]{amsart}
\usepackage{epsfig,graphics,amsmath,amssymb,amsthm}
\newtheorem{Theorem}{Theorem}[section]

\newtheorem{Corollary}[Theorem]{Corollary}
\newtheorem{Proposition}[Theorem]{Proposition}
\newtheorem{Lemma}[Theorem]{Lemma}

\theoremstyle{definition}
\newtheorem{Definition}[Theorem]{Definition}

\newtheorem{Remark}[Theorem]{Remark}

\begin{document}

\title[Automorphisms of complexes of curves]
{Automorphisms of complexes of curves on odd genus nonorientable
surfaces}

\author{Ferihe Atalan-Ozan}
\address{Department of Mathematics, Atilim University,
06836 \newline Ankara, TURKEY} \email{fatalan@atilim.edu.tr}
\date{\today}
\pagenumbering{arabic}

\begin{abstract}Let $N$ be a connected nonorientable surface of
genus $g$ with $n$ punctures. Suppose that $g$ is odd and $g+n
\geqslant 6$. We prove that the automorphism group of the complex
of curves of $N$ is isomorphic to the mapping class group
$\mathcal M_{N}$ of $N$.

\end{abstract}

\maketitle
\section{Introduction and statement of results}

Let $N$ be a connected nonorientable surface of genus $g$ with $n$
punctures and let $\mathcal M_N$ denote the mapping class group of
$N$, the group of isotopy classes of all diffeomorphisms $N\to N$.
The complex of curves $C(N)$ on $N$ is defined to be the abstract
simplicial complex whose vertices are the isotopy classes of
nontrivial simple closed curves.  A set of vertices
$\{v_{0},v_{1},\ldots,v_{q}\}$ forms a $q-simplex$ if and only if
$v_{0},v_{1},\ldots,v_{q}$ have pairwise disjoint representatives.

Each diffeomorphism $N \rightarrow N$ acts on the set of
nontrivial simple closed curves preserving the disjointness of
simple closed curves. It follows that the mapping class group of
$N$ acts on $C(N)$ as simplicial automorphisms. In other words,
there is a natural group homomorphism $\mathcal M_N \rightarrow
\textmd{Aut} \,C(N)$. The purpose of this paper is to prove the
following theorem.

\begin{Theorem}\label{Theorem-main1}
Let $N$ be a connected nonorientable surface of genus $g$ with $n$
punctures. Suppose that $g$ is odd and $g+n\geqslant 6$. Then the
natural map $\mathcal M_N \rightarrow {\rm Aut} \,C(N)$ is an
isomorphism.
\end{Theorem}

The complex of curves on an orientable surface $S$ was introduced
by Harvey~\cite{H}. It was shown by Ivanov~\cite{I2} and
Korkmaz~\cite{K1} that all automorphisms of the complex of curves
on $S$ are induced by diffeomorphisms of the surface $S$, with a
few exception. Another proof of this result was also obtained by
Luo ~\cite{L}. As a consequence of this result, it was proved that
any isomorphism between two subgroups of finite index in the
mapping class group of $S$ is given by the conjugation with a
mapping class (cf.~\cite{I2},~\cite{K1}). Ivanov also gave another
proof of the fact that the isometries of the Teichm\"uller space
are induced by diffeomorphisms of $S$.

Schaller~\cite{SS} showed that the extended mapping class group of
a hyperbolic surface is isomorphic to the automorphism group of
the graph. The set of vertices of the graph is the set of
nonseparating simple closed geodesic and the edges consisting of
pairs of those nonseparating curves satisfying the property that
the two curves in each pair intersect exactly once.
Margalit~\cite{M} proved that the automorphism group of the pants
complex is isomorphic to the extended mapping class group.
Irmak~\cite{Ir1,Ir2,Ir3} defined a superinjective simplicial map
and showed that a superinjective simplicial map of the complex of
curves and the complex of nonseparating curves is induced by a
homeomorphism of an orientable surface. Irmak and
Korkmaz~\cite{IK} proved that the automorphism group of the
Hatcher-Thurston complex is isomorphic to the extended mapping
class group modulo its center. Brendle and Margalit~\cite{BreM}
showed that any injection of a finite index subgroup of
$\mathcal{K}$, generated by the Dehn twists about separating
curves, into the Torelli group $\mathcal{I}$ of a closed
orientable surface is induced by a homeomorphism, conforming a
conjecture of Farb that $Comm(\mathcal{K}) \cong {\rm
Aut}(\mathcal{K}) \cong \mathcal M_S$. Behrstock and
Margalit~\cite{BehM} proved that for a torus with at least 3
punctures or a surface of genus 2 with at most 1 puncture, every
injection of a finite index subgroup of the extended mapping class
group into the extended mapping class group is the restriction of
an inner automorphism.

All the above results are about orientable surfaces. For
nonorientable surfaces we prove that the automorphism group of the
complex of curves of a nonorientable surface of odd genus is
isomorphic to its mapping class group. First, we show that the
natural group homomorphism is injective. Second, we prove that the
natural group homomorphism is surjective for the punctured
projective plane using the results and the ideas of Korkmaz's
analogous work on the punctured sphere \cite{K1} and a result by
Scharlemann \cite{S}. For higher genus, we use induction and some
ideas contained in Irmak's analogous work \cite{Ir2}.

\bigskip
\begin{Remark}\label{Remark-evengenus}
We believe that the same result holds for nonorientable surfaces
of even genus. Although we have some progress in proving the even
genus case, the proof is not complete yet.
\end{Remark}

\bigskip
\noindent{\bf Acknowledgement.}  This work is the author's
dissertation at Middle East Technical University. The author would
like to Mustafa Korkmaz, the author's thesis advisor, for his
continuous guidance, encouragement, suggestions and for reviewing
this manuscript.

\section{Preliminaries and notations}
\bigskip

Let $N$ be a connected nonorientable surface of genus $g$ with $n$
marked points. We call these marked points punctures. Recall that
the genus of a nonorientable surface is the maximum number of
projective planes in a connected sum decomposition.

\subsection{Circles and arcs}
If $a$ is a circle on $N$, by which we mean a simple closed curve,
then according to whether a regular neighborhood of $a$ is an
annulus or a Mobius strip, we call $a$ two-sided or one-sided
simple closed curve, respectively.

We say that a circle is nontrivial if it bounds neither a disc nor
annulus together with a boundary component, nor a disc with one
puncture, nor a Mobius band on $N$.

If $a$ is a circle, then we denote by $N_{a}$ the surface obtained
by cutting $N$ along $a$. A circle $a$ is called nonseparating if
$N_{a}$ is connected and separating otherwise. If $a$ is
separating, then $N_{a}$ has two connected components. If $a$ is
separating and if one of the components of $N_{a}$ is a disc with
$k$ punctures, then we say that $a$ is a $k$-separating circle.

We denote circles by the lowercase letters $a,b,c$ and their
isotopy classes by $\alpha,\beta,\gamma$. An embedded arc
connecting a puncture to itself or two different punctures will be
denoted by $a',b',c'$ and their isotopy classes by
$\alpha',\beta',\gamma'$.

Let $\alpha$ be the isotopy class of a circle $a$. We say that
$\alpha$ is nonseparating (respectively separating) if $a$ is
nonseparating (respectively separating). Similarly, we say that
$\alpha$ is one-sided, two-sided or $k$-separating vertex if $a$
is one-sided, two-sided or $k$-separating circle, respectively.

The geometric intersection number $i(\alpha,\beta)$ of two isotopy
classes $\alpha$ and $\beta$ is defined to be the infimum of the
cardinality of $a \cap b$ with $a \in \alpha$, $b \in \beta$. The
geometric intersection numbers $i(\alpha,\beta')$ and
$i(\alpha',\beta')$ are defined similarly.

The following lemma is proved in \cite{F}.
\begin{Lemma}\label{FLP}
Let $S$ be a sphere with $3$ punctures. Then
\begin{itemize}
\item[(i)]up to isotopy, there exists a unique nontrivial embedded
arc joining a puncture $P$ to itself, or $P$ to another puncture
$Q$, \item[(ii)]all circles on $S$ are trivial.
\end{itemize}

\end{Lemma}

\subsection{The complex of curves}
An abstract simplicial complex is defined as follows (cf.
~\cite{R}): Let $V$ be a nonempty set. An abstract simplicial
complex $K$ with vertices $V$ is a collection of nonempty finite
subsets of $V$, called simplices, such that if $v \in V$, then
$\{v\} \in K$, and if $\sigma \in K$ and $\sigma' \subset \sigma$
is a nonempty subset of $V$, then $\sigma' \in K$. The dimension
dim\,$\sigma$ of a simplex $\sigma$ is card\,$\sigma-1$, where
card\,$\sigma$ is the cardinality of $\sigma$. A simplex $\sigma$
is called a $q$-simplex if dim\,$\sigma = q$. The supremum of the
dimension of the simplices of $K$ is called the dimension of $K$.

A subcomplex $L$ of an abstract simplicial complex $K$ is called a
full subcomplex if whenever a set of vertices of $L$ is a simplex
in $K$, it is also a simplex in $L$.

The complex of curves $C(S)$ on an orientable surface $S$ is the
abstract simplicial complex whose vertices are the isotopy classes
of nontrivial simple closed curves. Similarly, the complex of
curves $C(N)$ on a nonorientable surface $N$ is the abstract
simplicial complex whose vertices are the isotopy classes of
nontrivial simple closed curves. In this complex of curves, we
take one-sided vertices as well as two-sided vertices. Clearly,
the complex of curves of a surface of genus $g$ with $n$ punctures
and with $b$ boundary components, and that of a surface of genus
$g$ with $n+b$ punctures are isomorphic. Therefore, sometimes we
regard boundary components and the punctures the same.

Two distinct vertices $\alpha , \beta \in C(N)$ are joined by an
edge if and only if their geometric intersection number is zero.
More generally, a set of vertices $\{v_0,v_1,\ldots,v_q \}$ forms
a $q$-simplex if and only if $i(v_j , v_k) = 0$ for all
$0\leqslant j,k \leqslant q $.

\subsubsection{Dimension} Clearly, the dimension of $C(N)$
is $n-2$ if $N$ is a projective plane with $n$ punctures. If $S$
is a sphere with $n$ punctures, then the dimension of $C(S)$ is
$n-4$. For higher genus, if $N$ is a connected nonorientable
surface of genus $g \geqslant 2$ with $n$ punctures such that the
Euler characteristic of $N$ is negative and if $g = 2r+1$, then
the dimension of $C(N)$ is $4r+n-2$ and if $g = 2r+2$, then the
dimension of $C(N)$ is $4r+n$ (see Section\,2.5). If $S$ is a
connected orientable surface of genus $g$ with $n$ punctures such
that $2g+n \geqslant 4$, then the dimension of $C(S)$ is $3g+n-4$.

\subsubsection{Links and dual links}
Let $\alpha$ be a vertex in the complex of curves. We define the
link $L(\alpha)$ of $\alpha$ to be the full subcomplex of the
complex of curves whose vertices are those of the complex of
curves which are joined to $\alpha$ by an edge in the complex of
curves. The dual link $L^{d}(\alpha)$ of $\alpha$ is the graph
whose vertices are those of $L(\alpha)$ such that two vertices of
$L^{d}(\alpha)$ are joined by edge if and only if they are not
joined by an edge in the complex of curves (or in $L(\alpha)$).

\subsubsection{Pentagons}
A pentagon is an ordered five-tuple $(\gamma_{1}, \gamma_{2},
\gamma_{3}, \gamma_{4}, \gamma_{5})$, defined up to cyclic
permutations and inversion, of vertices of the complex of curves
such that $i(\gamma_{i}, \gamma_{i+1}) = 0$ for $i = 1,2,\ldots,5$
and $i(\gamma_{i}, \gamma_{j}) \neq 0$ otherwise.

\subsection{Curve complexes of low dimensions} 
\begin{figure}[hbt]
\begin{center}
\includegraphics[width=10cm]{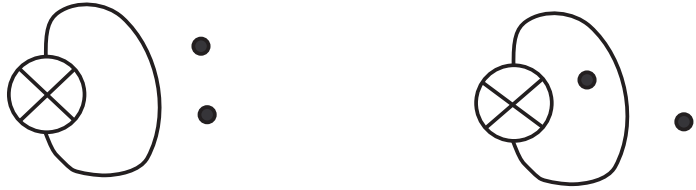}
\caption {The isotopy classes of circles on twice-punctured
$\mathbb{R}P^{2}$} \label{2puncture}
\end{center}
\end{figure}

Obviously, if $S$ is a
sphere with $\leqslant 3$, then there are no nontrivial circles on
$S$. Therefore, $C(S)$ is empty. If $S$ is a sphere with four
punctures, then $C(S)$ is infinite discrete. If $N$ is a
projective plane or a projective plane with one puncture then
$C(N)$ consists of a unique vertex. If $N$ is a projective plane
with two punctures, then $C(N)$ is finite (cf. \cite{S}). It
consists of two vertices, the isotopy classes of the circles
$c_{1}$ and $c_{2}$ of the Figure\,\ref{2puncture}.

\subsection{The arc complex $B(N)$} We now define another
abstract simplicial complex $B(N)$ as follows. The vertices of
$B(N)$ are the isotopy classes of nontrivial embedded arcs on $N$
connecting two punctures. A set of vertices of $B(N)$ forms a
simplex if and only if these vertices have pairwise disjoint
representatives.

\subsubsection{Arcs and $2$-separating circles}
If $a$ is $2$-separating circle on $N$, there exists up to isotopy
a unique nontrivial embedded arc $a'$ on the twice-punctured disc
component of $N_{a}$ joining two punctures by Lemma\,1.1 in
\cite{K1} and in \cite{F}. On the other hand, an arc $a'$
connecting two different punctures of $N$ determines uniquely a
$2$-separating circle up to isotopy, that is, the boundary of a
regular neighborhood of the arc. So, we have a one-to-one
correspondence between the set of $2$-separating isotopy classes
and the set of isotopy classes of embedded arcs connecting two
different punctures.

\subsubsection{Example}\label{Example-arcs}
\begin{figure}[hbt]
\begin{center}
\includegraphics[width=12cm]{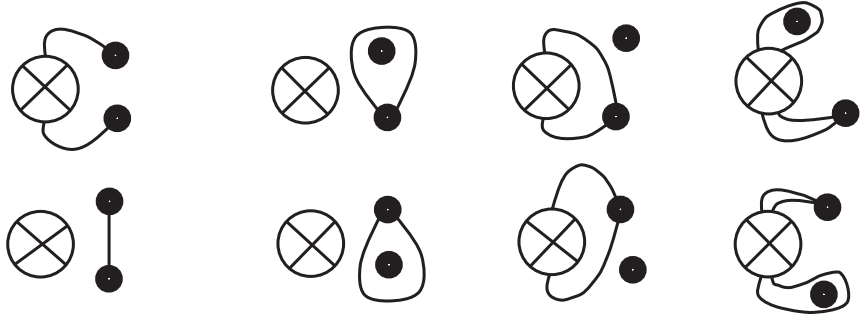}
\caption {The isotopy classes of arcs on twice-punctured
$\mathbb{R}P^{2}$} \label{2puncarc}
\end{center}
\end{figure}

If $N$ is a projective plane with one puncture, then $B(N)$
consists of a unique vertex. Scharlemann studied the arc complex
of a twice-punctured projective plane. He showed that the vertices
of the arc complex consist of the isotopy classes of arcs shown in
the Figure\,\ref{2puncarc}.

\subsubsection{Simple pairs and chains}
\begin{figure}[hbt]
\begin{center}
\includegraphics[width=12cm]{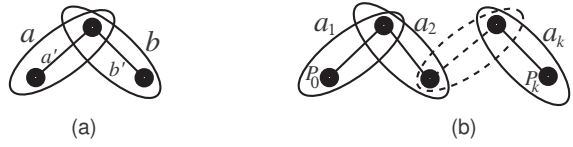}
\caption {A simple pair and a chain} \label{chain}
\end{center}
\end{figure}

If $a$ and $b$ are two $2$-separating circles, and $\alpha$,
$\beta$ are their isotopy classes, such that the corresponding
arcs $a'$ and $b'$ can be chosen disjoint with exactly one common
endpoint $P$, then we say that $a$ and $b$ constitute a simple
pair of circles and denote it by $\langle a ; b \rangle$ (see
Figure\,\ref{chain}(a)). Similarly, we say that $\langle a' ; b'
\rangle$ is a simple pair of arcs. We also say $\langle \alpha ;
\beta \rangle$ and $\langle \alpha' ; \beta' \rangle$ simple
pairs.

Let $a_{1}', a_{2}',\ldots, a_{k}'$ be embedded pairwise disjoint
arcs, $P_{i}$ and $P_{i+1}$ the endpoints of $a_{i}'$, with $P_{i}
\neq P_{j}$ for $i \neq j$, $0 \leqslant i, j \leqslant k$.
Therefore, $\langle a_{i}' ; a_{i+1}' \rangle$ is a simple pair of
arcs for each $1 \leqslant i \leqslant k-1$. Let $a_{1},
a_{2},\ldots, a_{k}$ be the corresponding circles. We say that
$\langle a_{1}'; a_{2}';\ldots; a_{k}' \rangle$ is a chain of arcs
(see Figure\,\ref{chain}(b)). Similarly, $\langle a_{1};
a_{2};\ldots; a_{k} \rangle$ is a chain of circles.

\subsection{Maximal simplices in the curve complex } We recall that
the maximum number of disjoint pairwise nonisotopic nontrivial
circles on a connected orientable surface $S$ of genus $r$ with
$n$ boundary components is $3r-3+n$, whenever the Euler
characteristic $\chi(S)$ of $S$ is negative.

\begin{Lemma}\label{Prop-odim1}
Let $S$ be a connected orientable surface of genus $g$ with $n$
punctures. Suppose that $2g+n \geqslant 4$. Then all maximal
simplices in $C(S)$ have the same dimension $3g+n-4$.
\end{Lemma}

\begin{Lemma}\label{Prop-odim2}
Let $N$ be a real projective plane with $n\geqslant 2$ punctures.
All maximal simplices in $C(N)$ have the same dimension $n-2$.
\end{Lemma}

\begin{proof}
Let $n = 2$. The complex $C(N)$ consists of only two vertices,
hence all simplices are of dimension $0$.

Let $n\geqslant 3$. We consider a maximal simplex $\sigma$ of
dimension $\ell$. Hence, $\sigma$ contains $\ell+1$ elements only
one of which is a one-sided vertex. By cutting $N$ along this
one-sided simple closed curve, we get sphere with $n+1$ punctures.
By Lemma~\ref{Prop-odim1}, all maximal simplices in the complex of
curves on the sphere with $n+1$ punctures have the same dimension
$n-3$. It follows that $\ell = n-2$.
\end{proof}

\begin{Proposition}\label{Prop-ndim1}
Let $N$ be a connected nonorientable surface of genus $g\geqslant
3$ with $n$ punctures such that the Euler characteristic of $N$ is
negative. Let $a_{r}=3r+n-2$ and $b_r=4r+n-2$ if $g=2r+1$, and
$a_{r}=3r+n-4$ and $b_r=4r+n-4$ if $g=2r$. Then there is a maximal
simplex of dimension $\ell$ in $C(N)$ if and only if $a_{r}
\leqslant \ell \leqslant b_{r}$.
\end{Proposition}

\begin{proof}
For each integer $\ell$ satisfying $a_{r} \leqslant \ell \leqslant
b_{r}$, the maximal simplicies are shown in
Figure\,\ref{maxsimplex0} and Figure\,\ref{maxsimplex1} for closed
nonorientable surface of odd genus. One can draw similar figures
for nonorientable surface of odd genus with punctures and
nonorientable surface of even genus. Moreover,
Figure\,\ref{maxsimplex} helps to see the maximal simplex of
dimension $\ell$ between $a_{r}$ and $b_{r}$.

\begin{figure}[hbt]
\begin{center}
\includegraphics[width=12cm]{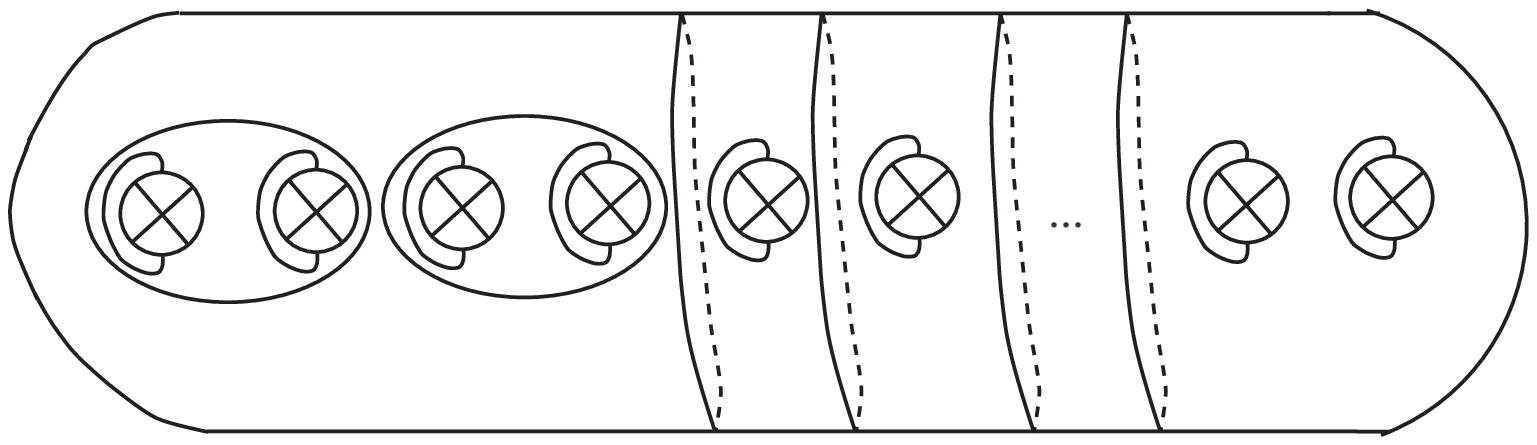}
\caption {maximum number of disjoint pairwise nonisotopic
nontrivial circles } \label{maxsimplex0}
\end{center}
\end{figure}

\begin{figure}[hbt]
\begin{center}
\includegraphics[width=12cm]{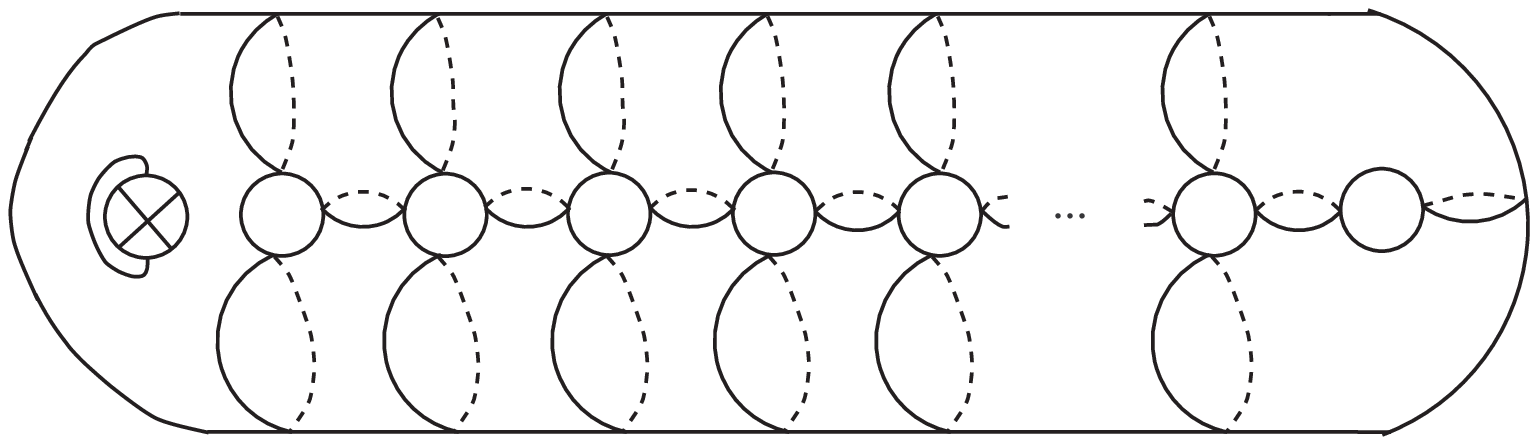}
\caption {maximum number of disjoint pairwise nonisotopic
nontrivial circles } \label{maxsimplex1}
\end{center}
\end{figure}

\begin{figure}[hbt]
\begin{center}
\includegraphics[width=12cm]{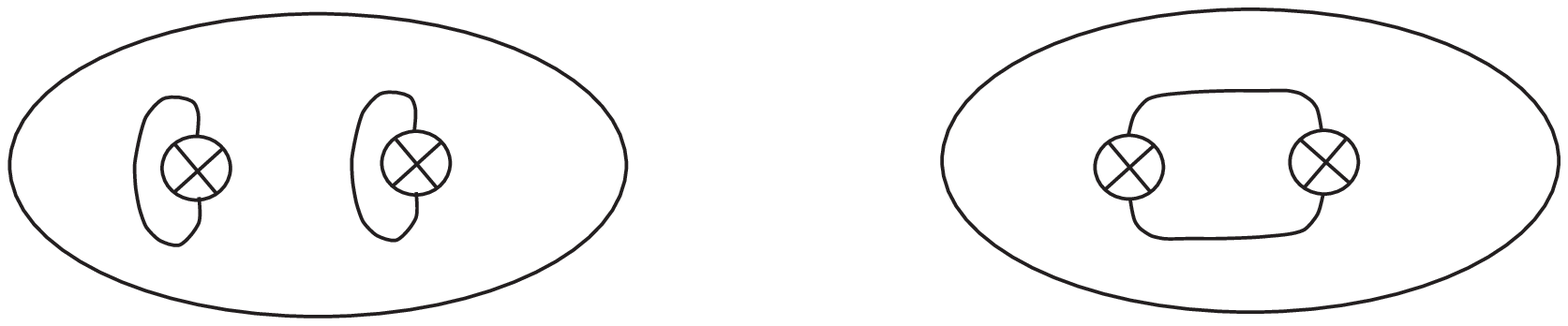}
\caption {} \label{maxsimplex}
\end{center}
\end{figure}

We now prove the converse. Let us consider a maximal simplex
$\sigma$ of dimension $\ell$. Hence, $\sigma$ contains $\ell +1$
elements. Choose pairwise disjoint simple closed curves
representing elements of $\sigma$, and let $N_\sigma$ denote the
surface obtained by cutting $N$ along these simple closed curves.

Suppose that the number of one-sided simple closed curves is $m$,
so that we have $\ell +1-m$ two-sided elements in $\sigma$. The
surface $N_\sigma$ is a disjoint union of $k$ pair of pants for
some positive integer $k$. By the Euler characteristic argument,
it can be seen that $k=g+n-2$. The number of boundary components
and punctures on $N_\sigma$ is $3k$. By counting the contribution
of one-sided curves and two-sided curves to the boundary of
$N_\sigma$, one can easily see that
 \begin{eqnarray} \label{eqn3k}
 3k=n+m+2(\ell+1-m).
 \end{eqnarray}

Suppose first that $g=2r+1$. In this case $1\leq m\leq 2r+1$. From
the equality~(\ref{eqn3k}), it is easy to see that $m$ must be odd
and $\ell=3r+n-2 + {{m-1}\over 2}$.

Suppose now that $g=2r$. In this case $0\leq m\leq 2r$. From the
equality~(\ref{eqn3k}), it is easy to see that $m$ must be even
and $\ell=3r+n-4 + {{m}\over 2}$.

The proposition follows from these.
\end{proof}

\subsection{Centralizer of certain subgroups}

\begin{figure}[hbt]
\begin{center}
\includegraphics[width=12cm]{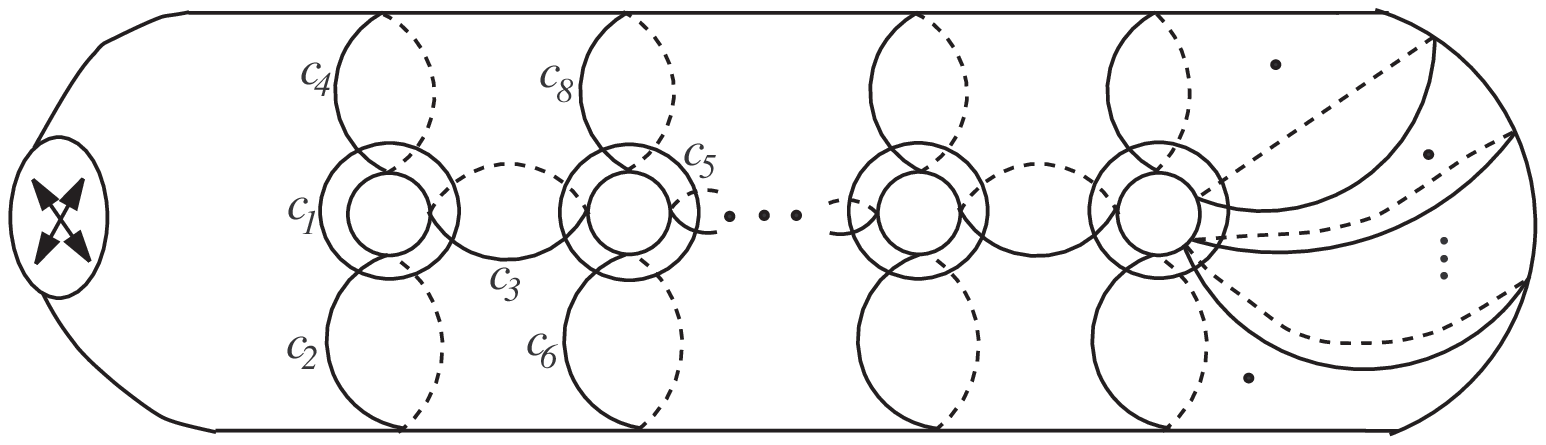}
\caption {} \label{surf0-circles}
\end{center}
\end{figure}

Let $\mathcal{T'}$ be the subgroup of mapping class group of $N$
such that $\mathcal{T'}$ is generated by the Dehn twist of
two-sided nonseperating circles as below shown in
Figure\,\ref{surf0-circles}.

\begin{Proposition}\label{Prop-centralizer}
Let $N$ be a connected nonorientable surface of odd genus $g$ and
$g+n \geqslant 5$. Suppose that $C$ is a collection of two-sided
nonseparating circles in Figure\,\ref{surf0-circles} such that
$\mathcal{T'}$ is generated by the Dehn twist $t_{c_{i}}$ along
the circles $c_{i}$ of $C$. Then the centralizer
$\mathcal{C}_{\mathcal M_{N}}(\mathcal{T'})$ is trivial.
\end{Proposition}

\begin{proof}
Let $[f] \in \mathcal{C}_{\mathcal M_{N}}(\mathcal{T'})$.
Unoriented circles will be denoted by $\bar{c}$. Then,
$f(\bar{c_{i}})=\bar{c_{i}}$ for $c_{i} \in C$. Let $a$ be a
one-sided circle such that $N_{a}$ is an orientable surface. We
observe that $f(\bar{a})=\bar{a}$. Because there is one isotopy
class of one-sided circle which does not intersect circles $c_{i}$
in $C$. Therefore, by cutting $N$ along $a$, we get diffeomorphism
$f_{|} : N_{a} \rightarrow N_{a}$. Moreover, $f_{|}$ is
orientation preserving. To see this, assume that $c_{j}$ is dual
to $c_{i}$, we know that $f_{|}$ preserves the orientation of a
regular neighborhood of $c_{i}$ if and only if it preserves the
orientation of $c_{j}$. Recall that tubular neighborhood of $c_{i}
\cup c_{j}$ is a torus with one boundary component. Since the
product of the orientations of $c_{i}$ and $c_{j}$ gives the
orientation of the torus, $f_{|}$ preserves the orientation of
tubular neighborhood. Using this argument for these circles in
Figure\,\ref{surf0-circles}, we deduce that $f_{|}$ preserves the
orientation of the surface $N_{a}$. Since $\mathcal{C}_{\mathcal
M_{N_{a}}}(\mathcal PM_{N_{a}})=\{1\}$ in \cite{I3}, $f_{|}$ is
isotopic to identity. Since $f(a)=a$, we obtain that $f$ is
isotopic to identity on $N$. Hence, $\mathcal{C}_{\mathcal
M_{N}}(\mathcal{T'})=\{1\}$.
\end{proof}

\bigskip
\section{Injectivity of $\mathcal M_N \rightarrow \textmd{Aut}
\,C(N)$}

We first show that the natural map $\mathcal M_N \rightarrow
\textmd{Aut} \,C(N)$ is injective.

\begin{Theorem}\label{Lemma-injec0}
Let $N$ be a connected nonorientable surface of genus $g$ with $n$
punctures. Suppose that $g+n \geqslant 6$. Then the natural map
$\mathcal M_N \rightarrow {\rm Aut} \,C(N)$ is injective.
\end{Theorem}

\begin{proof}
Suppose that $[f]\in \mathcal M_N$ acts trivially on $C(N)$. Then
$f : N \rightarrow N$, $f_{*} : C(N) \rightarrow C(N)$, $f_{*}=
id$; i.e., $f_{*}(\nu) = \nu $ for all $\nu \in C(N)$. First,
suppose that $g=2r+1$ is odd. Let $a$ be a one-sided circle such
that $N_{a}$ is orientable surface of genus $r$ with $n+1$
boundary components. We denote by $\alpha$ the isotopy class of
$a$. Then $f_{*}(\alpha)= \alpha$. This implies $f(a)$ is isotopic
to $a$. Then there exists a diffeomorphism $g \backsimeq id$,
$g(f(a)) = a$. Let $h = g \circ f$. We observe that $h(a) = a$ and
$h_{*} = g_{*} \circ f_{*} = id \circ id = id$. Therefore, we have
a diffeomorphism $h_{|} : N_{a} \longrightarrow N_{a}$ such that
$(h_{|})_{*} = id$. Using Ivanov's Theorem in \cite{I2} (for $r=0$
and $1$, using Theorem\,1 in \cite{K1}), we get $h_{|} \backsimeq
id$. Since $h(a) = a$, we see that $h \backsimeq id$ on $N$. In
other words, it descends to a diffeomorphism of $N$. So, we have
$g \circ f \backsimeq id$. Since $g \backsimeq id$, we get $f
\backsimeq id$. Hence, the natural map $\mathcal M_N \rightarrow
\textmd{Aut} \,C(N)$ is injective.

Now, suppose that $g=2r+2$. Let $b$ be a nonseparating two-sided
circle such that $N_{b}$ is orientable surface of genus $r$ with
$n+2$ boundary components. Let $\beta$ be the isotopy class of
$b$. Then $f_{*}(\beta)= \beta$. This implies $f(b)$ is isotopic
to $b$. Then there exists a diffeomorphism $g \backsimeq id$,
$g(f(b)) = b$. Let $h = g \circ f$. We see that $h(b) = b$ and
$h_{*} = g_{*} \circ f_{*} = id \circ id = id$. Therefore, we have
a diffeomorphism $h_{|} : N_{b} \longrightarrow N_{b}$ such that
$(h_{|})_{*} = id$. Using Ivanov's Theorem in \cite{I2} (for $r=0$
and $1$, using Theorem\,1 in \cite{K1}), we get $h_{|} \backsimeq
id$. Since $h(b) = b$, $h \backsimeq id$ on $N$. In other words,
it descends to a diffeomorphism of $N$. Therefore, we obtain that
$g \circ f \backsimeq id$. Since $g \backsimeq id$, we get $f
\backsimeq id$. Hence, the natural map $\mathcal M_N \rightarrow
\textmd{Aut} \,C(N)$ is injective.
\end{proof}

\section{Punctured $\mathbb{R}P^{2}$}
\bigskip

Throughout this section unless otherwise stated, $N$ will denote a
real projective plane with $n \geqslant 5$ punctures. We need at
least $5$ punctures for the proof of the Lemma~\ref{Lemma-indep}.

In this section, we will prove that the natural homomorphism
$\mathcal M_N \rightarrow \textmd{Aut}\,C(N)$ is surjective. Hence
it will be an isomorphism. For this, we first prove that
automorphisms of $C(N)$ preserve the topological type of the
vertices of $C(N)$ and that certain pairs of vertices of $C(N)$
can be realized in the complex of curves. We conclude that
automorphisms of $C(N)$ preserve these pairs of vertices. Next, we
show that every automorphism of $C(N)$ induces an automorphism of
the complex $B(N)$ in a natural way. The automorphisms of $B(N)$
are determined by their action on a maximal simplex. Then, we use
the relation between maximal simplices of $B(N)$ and isotopy
classes of ideal triangulations of $N$ and in conclusion, we show
that an automorphism of $B(N)$ induced by some automorphism of
$C(N)$ agrees with a mapping class.

We remind that if $N$ is a projective plane with $n \geqslant 2$
punctures, then up to diffeomorphism there is only one
nonseparating one-sided circle and also there is no nonseparating
two-sided circle. The other circles are $k$-separating for some
$k$.
\bigskip

\begin{Lemma}\label{Lemma-iso}
Let $n \geqslant 2$ and $k\geqslant 4$. If $N$ is a projective
plane with $n$ punctures and $S$ is a sphere with $k$ punctures,
then $C(N)$ and $C(S)$ are not isomorphic.
\end{Lemma}

\begin{proof}
The complexes $C(N)$ and $C(S)$ have dimensions $n-2$ and $k-4$,
respectively. If $k \neq n+2$, since these complexes of curves
have different dimensions, $C(N)$ and $C(S)$ are not isomorphic.
If $k = n+2$, we proceed as follows.

Let $n=2$. Then $C(N)$ is finite (see \cite{S}), however, $C(S)$
is infinite discrete since $S$ is a sphere with $4$ punctures.
Therefore, they are not isomorphic.

Now, assume that $n \geqslant 3$ and $C(N)$ and $C(S)$ are not
isomorphic when $N$ has $n-1$ punctures. We need to show that
these complexes are not isomorphic if $N$ has $n$ punctures.
Assume that there is an isomorphism $\varphi : C(N) \rightarrow
C(S)$. Note that for a vertex $\gamma$ of $C(N)$, the dual link
$L^{d}(\gamma)$ of $\gamma$ is connected if and only if $\gamma$
is either one-sided or $2$-separating. For a vertex $\delta$ of
$C(S)$, the dual link $L^{d}(\delta)$ of $\delta$ is connected if
and only if $\delta$ is $2$-separating. From this, it follows that
the image of the union of the set of one-sided vertices and the
set of $2$-separating vertices of $C(N)$ is precisely the set of
$2$-separating vertices of $C(S)$. Let $\gamma$ be a
$2$-separating vertex of $C(N)$. Then $\varphi$ takes $\gamma$ to
a $2$-separating vertex $\delta$ of $C(S)$ and induces an
isomorphism $L(\gamma) \rightarrow L(\delta)$.

Clearly, $L(\gamma)$ is isomorphic to the complex of curves of a
real projective plane with $n-1$ punctures and $L(\delta)$ is
isomorphic to the complex of curves of a sphere with $n+1$
punctures. By assumption, these complexes are not isomorphic. We
get a contradiction. Hence, $C(N)$ and $C(S)$ are not isomorphic.
\end{proof}

\begin{Theorem}\label{Theorem-prese}
The group ${\rm Aut} \,C(N)$ preserves the topological type of the
vertices of $C(N)$.
\end{Theorem}

\begin{proof}
Let $\varphi$ be an automorphism of $C(N)$. Note that for a vertex
$\gamma$ of $C(N)$, the dual link $L^{d}(\gamma)$ of $\gamma$ is
connected if and only if $\gamma$ is either one-sided vertex or
$2$-separating. Therefore, $\varphi$ cannot take a one-sided
vertex or a $2$-separating to a $k$-separating vertex with $k >
2$.

Assume that $\varphi(\alpha)=\beta$ is a $2$-separating vertex for
some one-sided vertex $\alpha$. Let $a \in \alpha$ be  a circle.
Then $N_{a}$ is disc with $n$ punctures. Let $b \in \beta$. Then
$N_{b}$ is homeomorphic to $\mathbb{R}P^{2}$ with $n-2$ punctures
and with one boundary component. Clearly, $L(\alpha)$ is
isomorphic to $C(N_{a})$. Similarly,  $L(\beta)$ is isomorphic to
$C(N_{b})$. Since the complexes $C(N_{a})$ and $C(N_{b})$ are not
isomorphic by Lemma~\ref{Lemma-iso}, $\beta$ cannot be
$2$-separating. It also follows from this that $\varphi$ maps a
$2$-separating vertex to a $2$-separating vertex.

Let $\gamma$ be a $k$-separating vertex for some $3 \leqslant k
\leqslant n-1$. Then $\varphi(\gamma)=\delta$ is an $l$-separating
vertex for some $3 \leqslant l \leqslant n-1$. We must show that
$k = l$. Let $c \in \gamma$ and $d \in \delta$ be circles. $N_{c}$
and $N_{d}$ have two connected components. Let $N_{c} = N_{0}
\sqcup N_{1}$ and $N_{d} = N'_{0} \sqcup N'_{1}$ such that $N_{i}$
and $N'_{i}$ have genera $i$. The dual link $L^{d}(\gamma)$ has
exactly two connected components. Let us denote these components
by $L_{0}^{d}(\gamma)$ and $L_{1}^{d}(\gamma)$; $L^{d}(\gamma) =
L_{0}^{d}(\gamma) \sqcup L_{1}^{d}(\gamma)$. We name
$L_{i}(\gamma)$ so that the vertices of $L^{d}_{i}(\gamma)$ are
the isotopy classes of circles on $N_{i}$. Let $L_{i}(\gamma)$ be
the full subcomplex of $C(N)$ with vertices $L^{d}_{i}(\gamma)$.
It follows that $(L_{i}(\gamma))^{d} = L^{d}_{i}(\gamma)$. Then
$L_{i}(\gamma)$ is isomorphic to $C(N_{i})$. Clearly, the
dimension of $C(N_{0})$ is $k-3$ since $N_{0}$ is a sphere with
punctures and the dimension of $C(N_{1})$ is $n-k-2$. Similarly,
we define $L^{d}_{i}(\delta)$ and $L_{i}(\delta)$ so that
$L_{i}(\delta)$ is isomorphic to $C(N'_{i})$. If $\varphi(\gamma)
= \delta$, then $\varphi$ restricts to an isomorphism $ L(\gamma)
\rightarrow L(\delta)$, which induces an isomorphism $
L^{d}(\gamma) \rightarrow L^{d}(\delta)$. Then either
$\varphi(L_{0}^{d}(\gamma)) = L_{0}^{d}(\delta)$ or
$\varphi(L_{0}^{d}(\gamma)) = L_{1}^{d}(\delta)$. Hence,
$\varphi(L_{0}(\gamma)) = L_{0}(\delta)$ or
$\varphi(L_{0}(\gamma)) = L_{1}(\delta)$. However,
$\varphi(L_{0}(\gamma)) = L_{1}(\delta)$ is not possible by
Lemma~\ref{Lemma-iso}. Therefore, $\varphi(L_{0}(\gamma)) =
L_{0}(\delta)$. It follows that their dimensions are equal: $k-3 =
l-3$. Hence, $k = l$.

The proof of the theorem is now complete.
\end{proof}

\begin{Theorem}[Korkmaz]\label{Korkmaz-penta}
Let $\alpha$ and $\beta$ be two $2$-separating vertices of the
complex of curves $C(S)$. Then $\langle\alpha;\beta\rangle$ is a
simple pair if and only if there exist vertices
$\gamma_{1},\gamma_{2},\gamma_{3},\ldots,\gamma_{n-2}$ of $C(S)$
satisfying the following conditions.

\begin{itemize}

\item[(i)] $(\gamma_{1},\gamma_{2},\alpha,\gamma_{3},\beta)$ is a
pentagon in $C(S)$,

\item[(ii)] $\gamma_{1}$ and $\gamma_{n-2}$ are $2$-separating,
$\gamma_{2}$ is $3$-separating, and $\gamma_{k}$ and
$\gamma_{n-k}$ are $k$-separating for $3 \leqslant k \leqslant
n/2$,

\item[(iii)] $\{\alpha,\gamma_{3}\}\cup\sigma$,
$\{\alpha,\gamma_{2}\}\cup\sigma$,
$\{\beta,\gamma_{3}\}\cup\sigma$ and
$\{\gamma_{1},\gamma_{2}\}\cup\sigma$ are codimension-zero
simplices, where
$\sigma=\{\gamma_{4},\gamma_{5},\ldots,\gamma_{n-2} \}$.

\end{itemize}
\end{Theorem}

\begin{Theorem}\label{Theorem-penta}
Let $\alpha$ and $\beta$ be two $2$-separating vertices of $C(N)$.
Then $\langle\alpha;\beta\rangle$ is a simple pair if and only if
there exist vertices
$\gamma_{1},\gamma_{2},\gamma_{3},\ldots,\gamma_{n-1},\delta$ of
$C(N)$ satisfying the following conditions.

\begin{itemize}

\item[(i)] $(\gamma_{1},\gamma_{2},\alpha,\gamma_{3},\beta)$ is a
pentagon in $C(N)$,

\item[(ii)] $\gamma_{1}$ is $2$-separating, $\gamma_{2}$ is
$3$-separating, and $\gamma_{k}$ is $k$-separating for $3
\leqslant k \leqslant n-1$, $\delta$ is one-sided,

\item[(iii)] $\{\alpha,\gamma_{3}\}\cup\sigma \cup \{ \delta \}$,
$\{\alpha,\gamma_{2}\}\cup\sigma \cup\{ \delta \}$,
$\{\beta,\gamma_{3}\}\cup\sigma \cup\{ \delta \}$ and
$\{\gamma_{1},\gamma_{2}\}\cup\sigma \cup\{ \delta \}$ are
codimension-zero simplices, where
$\sigma=\{\gamma_{4},\gamma_{5},\ldots,\gamma_{n-1} \}$.

\end{itemize}
\end{Theorem}

\begin{proof}
Suppose that $\langle \alpha ; \beta \rangle$ is a simple pair.
Let $a \in \alpha$ and $b \in \beta$ so that $\langle a ; b
\rangle$ is a simple pair. It is clear that any two simple pairs
of circles are topologically equivalent, that is; if $\langle c ;
d \rangle$ is any other simple pair, then there exists a
diffeomorphism $F:N \rightarrow N$ such that $\langle F(c) ; F(d)
\rangle = \langle a ; b \rangle$. So, we can assume that $a$ and
$b$ are the circles illustrated in Figure\,\ref{penta}. In the
figure, we think of the sphere as the one point compactification
of the plane and the cross inside the circle means that we delete
open disc and identify the antipodal boundary points so that we
get a real projective plane with punctures. The isotopy classes
$\gamma_{i}$ of the circles $c_{i}$ and the isotopy class $\delta$
of the one-sided circle $d$ satisfy (i)-(iii).

\begin{center}
\begin{figure}[hbt]
\includegraphics[width=12cm]{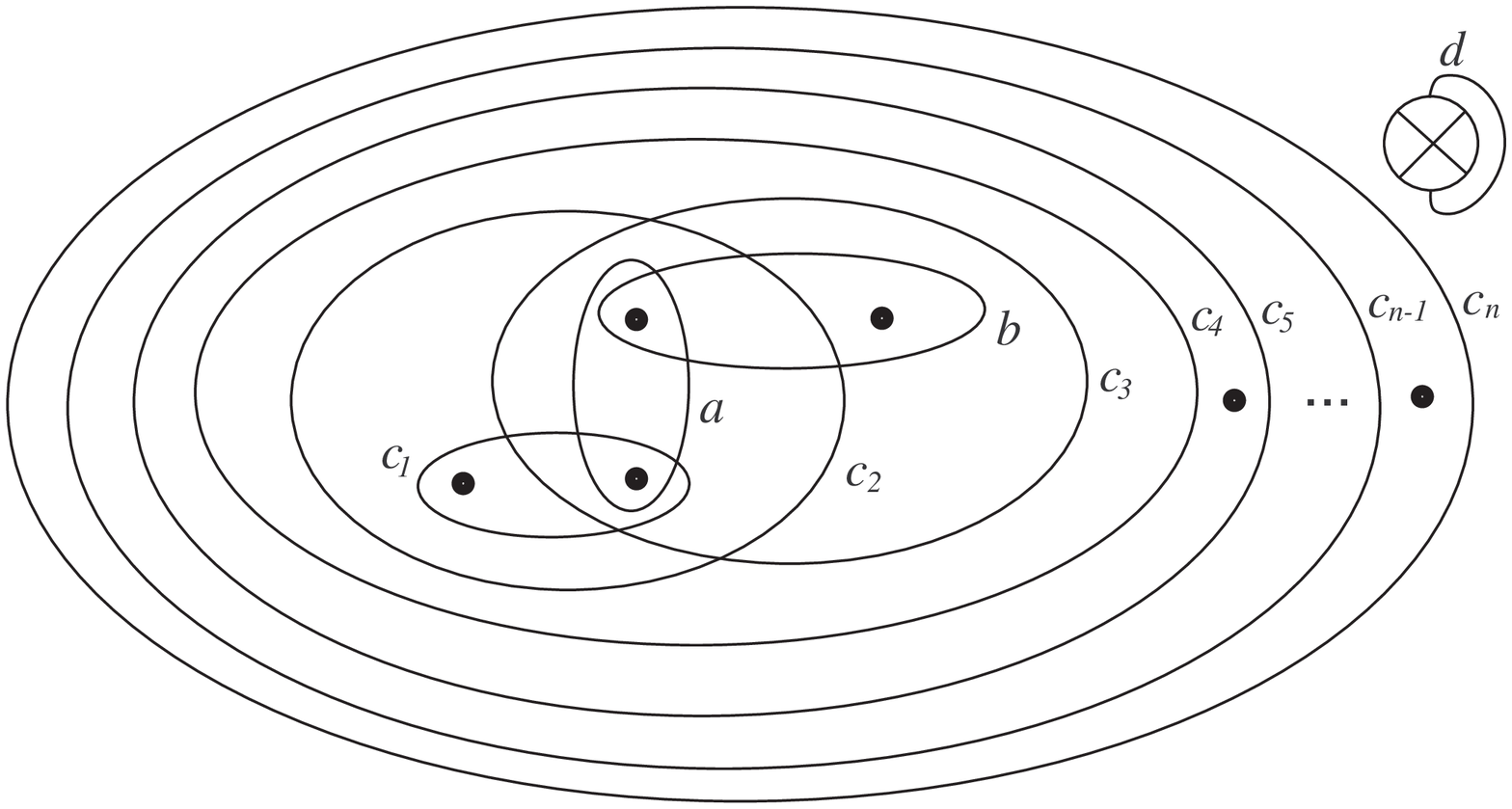}
\caption{} \label{penta}
\end{figure}
\end{center}

Now, we prove the converse. Assume that conditions (i)-(iii) above
hold. Let $d \in \delta$ be a one-sided circle. Deleting
$\{\delta\}$ from conditions (ii) and (iii), we have
codimension-one simplices $\{\alpha,\gamma_{3}\}\cup\sigma$,
$\{\alpha,\gamma_{2}\}\cup\sigma$,
$\{\beta,\gamma_{3}\}\cup\sigma$ and
$\{\gamma_{1},\gamma_{2}\}\cup\sigma$. However, these simplicies
are codimension-zero simplices in the complex $C(N_{d})$. By
Theorem\,\ref{Korkmaz-penta}, we see that
$\langle\alpha;\beta\rangle$ is a simple pair on the sphere and
$N_{d}$ is a sphere with $n+1 \geqslant 6$ punctures. Let us say
that a puncture is inside $a$ if it is one of the two punctures on
the disc bounded by $a$. Similarly for $b$. There are three
possibilities the boundary component $d'$ (we see it as a puncture
from point of view of the curve complex) resulting from cutting
along $d$, as illustrated in Figure\,\ref{pair}. In the figure,
this boundary component $d'$ is drawn as an oval. The case are:
$(1)$ it may be outside of both $a$ and $b$, $(2)$ it may be
inside, say, $a$ and outside of $b$, or $(3)$ it may be the unique
puncture inside both $a$ and $b$.

\begin{figure}[hbt]
\begin{center}
\includegraphics[width=12cm]{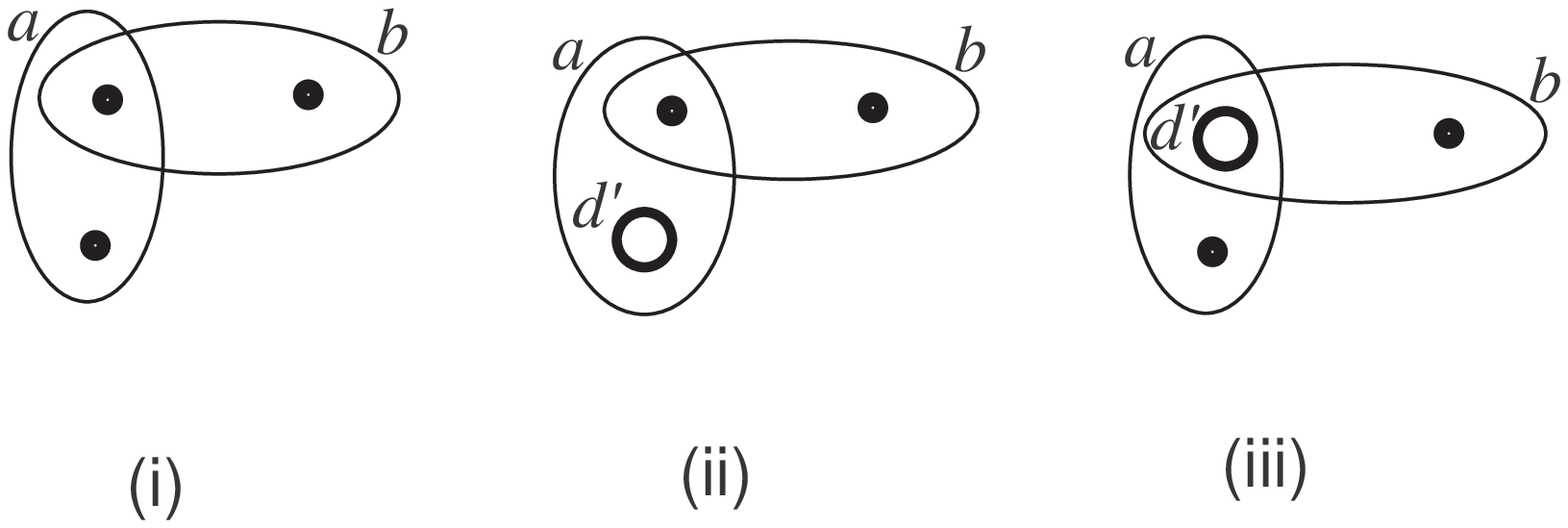}
\caption{} \label{pair}
\end{center}
\end{figure}

After the identifying the antipodal boundary points of the oval,
in order to get both $a$ and $b$ as $2$-separating circles, we
must have the case in Figure\,\ref{pair}(i). Hence, we see that
$\langle\alpha;\beta\rangle$ is a simple pair on $N$ as well.
\end{proof}

\begin{Corollary}\label{Corollary-chain}
Let $\varphi$ be an automorphism of $C(N)$. If
$\langle\alpha;\beta\rangle$ is a simple pair, so is $\langle
\varphi(\alpha) ; \varphi(\beta) \rangle$. Similarly, the image of
a chain in $C(N)$ under $\varphi$ is also a chain.
\end{Corollary}

\begin{proof}
The conditions (i) and (iii) of Theorem\,\ref{Theorem-penta} are
invariant under the automorphisms of $C(N)$. By
Theorem\,\ref{Theorem-prese}, the condition (ii) is also invariant
under the automorphisms of $C(N)$. Corollary follows from these.
\end{proof}

\bigskip
\subsubsection{Action of ${\rm Aut} \,C(N)$ on punctures}

We can define an action of the group $\textmd{Aut} \,C(N)$ on
punctures of $N$ as follows. For $\varphi \in \textmd{Aut} \,C(N)$
and for a puncture $P$ of $N$, choose two isotopy classes
$\alpha'$, $\beta'$ of embedded arcs such that $\langle
\alpha';\beta' \rangle$ is a simple pair with the common endpoint
$P$. Define $\varphi(P)$ to be the common endpoint of the simple
pair $\langle \varphi(\alpha'); \varphi(\beta') \rangle$. We note
that by the one-to-one correspondence between the set of
$2$-separating vertices of $C(N)$ and the set of those vertices of
$B(N)$ which join different punctures, $\textmd{Aut} \,C(N)$ has a
well-defined action on the latter set.

Next two lemmas below can be proved in the same way as Lemma\,3.5
and Lemma\,3.6 in \cite{K1}.

\begin{Lemma}\label{Lemma-indep}
The definition of the action of ${\rm Aut} \,C(N)$ on the
punctures of $N$ is independent of the choice of the simple pair.
\end{Lemma}

\begin{Lemma}\label{Lemma-punc}
Let $\varphi \in {\rm Aut} \,C(N)$, $\alpha$ be a $k$-separating
vertex of $C(N)$ and $a \in \alpha$. If $N_{a}'$ denotes the
$k$-punctured disc component of $N_{a}$ and $N_{a}''$ denotes
$(n-k)$-punctured $\mathbb{R}P^{2}$ with one boundary component,
then $\varphi(\mathcal{P}(N_{a}')) = \mathcal{P}(N_{\varphi(a)}')$
and $\varphi(\mathcal{P}(N_{a}'')) =
\mathcal{P}(N_{\varphi(a)}'')$.
\end{Lemma}

\bigskip
\subsubsection{Action of ${\rm Aut} \,C(N)$ on arcs}

We define  an action of $\textmd{Aut} \,C(N)$ on the vertices of
$B(N)$ as follows. Let $\varphi \in \textmd{Aut} \,C(N)$,
$\alpha'$ a vertex of $B(N)$ and let $a' \in \alpha'$. If $a'$ is
joining two different punctures, then $\varphi(\alpha')$ is
already defined by the correspondence between the $2$-separating
vertices of $C(N)$ and the action of $\textmd{Aut} \,C(N)$ on
$C(N)$. The other words, $\varphi(\alpha')$ is the isotopy class
of the arc, which is unique up to isotopy, joining two punctures
on the twice-punctured disc component of $N_{\varphi(a)}$ for
$\varphi(a) \in \varphi(\alpha)$.

Suppose now that the arc $a'$ joins a puncture $P$ to itself  and
it is a two-sided loop. Then the action is same as in \cite{K1}.
To be more precise, let $a_{1}$ and $a_{2}$ be the boundary
components of a regular neighborhood of $a' \cup \{P\}$ and
$\alpha_{1}$ and $\alpha_{2}$ be their isotopy classes,
respectively. Since $a'$ cannot be deformed to $P$, at most one of
$a_{1}$ and $a_{2}$ is trivial.

If $a_{1}$ is trivial, then $a_{1}$ bounds either a disc with one
puncture or a M\"obius band. If it bounds a disc with a puncture
$Q$, then $a_{2}$ bounds a disc with two punctures $P$ and $Q$. By
Theorem~\ref{Theorem-prese}, $\varphi(\alpha_{2})$ is
$2$-separating. So for a representative $\varphi(a_{2})$ of
$\varphi(\alpha_{2})$, one of the components, say
$N_{\varphi(a_{2})}'$ of $N_{\varphi(a_{2})}$ is a twice-punctured
disc, with punctures $\varphi(P)$ and $\varphi(Q)$ by
Lemma~\ref{Lemma-punc}. Define $\varphi(\alpha')$ to be the
isotopy class of a nontrivial simple arc on $N_{\varphi(a_{2})}'$
joining $\varphi(P)$ to itself. Such an arc is unique up to
isotopy by Lemma~\ref{FLP}.

If $a_{1}$ bounds a M\"obius band, then $a_{2}$ bounds a
projective plane with one puncture $P$. By
Theorem~\ref{Theorem-prese}, $\varphi(\alpha_{2})$ bounds a
projective plane with one puncture $\varphi(P)$. Therefore, for a
representative $\varphi(a_{2})$ of $\varphi(\alpha_{2})$, one of
components, say $N'_{\varphi(a_{2})}$, of $N_{\varphi(a_{2})}$ is
a projective plane with one puncture and one boundary component.
The puncture on $N'_{\varphi(a_{2})}$ is $\varphi(P)$ by
Lemma~\ref{Lemma-punc}. There is only one nontrivial two-sided
loop (arc) joining a puncture $\varphi(P)$ to itself by
Example~\ref{Example-arcs}. We define $\varphi(\alpha')$ to be the
isotopy class of this two-sided loop (arc) joining $\varphi(P)$ to
itself.

If neither of $a_1$ and $a_2$ is trivial, then $a_{1}$ and $a_{2}$
bound an annulus with a puncture $P$.  We claim that
$\varphi(a_{1})$ and $\varphi(a_{2})$ also bound a once-punctured
annulus with only one puncture $\varphi(P)$. Here,
$\varphi(a_{i})$ is a representative of $\varphi(\alpha_{i})$ for
$i=1,2$. To see this, let $N_{a_{i}}'$ be the subsurface of $N$
bounded by $a_i$ not containing the puncture $P$. Similarly we
define $N'_{\varphi ({a_{i}})}$ to be the component of $N_{\varphi
({a_{i}})}$ not containing $\varphi(P)$. Now, assume that the set
of punctures on $N_{a_{1}}'$ and $N_{a_{2}}'$ are
$\mathcal{P}(N_{a_{1}}') = \{P_{1},\ldots,P_{k}\}$ and
$\mathcal{P}(N_{a_{2}}') = \{P_{k+1},\ldots,P_{n-1}\}$,
respectively. Then $P_{i} \neq P_{j}$ for all $i,j$. By
Lemma~\ref{Lemma-punc}, $\mathcal{P}(N_{\varphi(a_{1})}') =
\{\varphi(P_{1}),\ldots,\varphi(P_{k})\}$ and
$\mathcal{P}(N_{\varphi(a_{2})}') =
\{\varphi(Q_{1}),\ldots,\varphi(Q_{n-k})\}$. We deduce that since
$\varphi(a_{1})$ and $\varphi(a_{2})$ are disjoint and
nonisotopic, they must bound an annulus with only one puncture
$\varphi(P)$. The class $\varphi(\alpha')$ is defined to be the
isotopy class of the unique arc up to isotopy on this annulus
joining $\varphi(P)$ to itself.

Suppose finally that $a'$ is a one-sided loop (arc) joining a
puncture $P$ to itself. Let $a$ be the boundary component of a
regular neighborhood of $a'\cup\{P\}$ and let $\alpha$ be the
isotopy class of $a$. The circle $a$ bounds a M\"obius band with a
puncture $P$. By Theorem~\ref{Theorem-prese}, $\varphi(a) \in
\varphi(\alpha)$ bounds a M\"obius band $M$ with a puncture, say
$Q$. By Lemma~\ref{Lemma-punc}, $Q = \varphi(P)$. By
Example~\ref{Example-arcs}, there is up to isotopy a unique
one-sided loop $b'$ on $M$ joining $\varphi(P)$ to itself. We
define $\varphi(\alpha')$ to be the isotopy class of $b'$.

\begin{Lemma}\label{Lemma-disjoint}
Let $\varphi$ be an automorphism of $C(N)$ and $\alpha'$ and
$\beta'$ be two distinct vertices of $B(N)$ such that $i(\alpha',
\beta') = 0$. Then $i(\varphi(\alpha'), \varphi(\beta')) = 0$.
Therefore, every automorphism of $C(N)$ yields an automorphism of
$B(N)$.
\end{Lemma}

\begin{proof}
Let $a'$ and $b'$ be two disjoint representatives of $\alpha'$ and
$\beta'$, respectively. There are thirteen possible cases as
illustrated in Figure\,\ref{puncture}. In each figure, we assume
that the arc on the left is $a'$ and the other is $b'$.

If $a'$ (respectively $b'$) is joining two different punctures, we
denote by $\alpha$ (respectively $\beta$) the $2$-separating
vertex of $C(N)$ corresponding to $\alpha'$ (respectively
$\beta'$), and by $a$ (respectively $b$) a representative of
$\alpha$(respectively $\beta$).

If $a'$ (respectively $b'$) is a two-sided arc  joining a puncture
$P$ to itself, we denote by $a_{1}$ and $a_{2}$ (respectively
$b_{1}$ and $b_{2}$) the boundary components of a regular
neighborhood of $a' \cup \{P\}$ (respectively $b' \cup \{P\}$). If
a regular neighborhood of $a' \cup \{P\}$ (respectively $b' \cup
\{P\}$) is a once-punctured M\"obius band $M$, we denote the
boundary component of $M$ by $a_{3}$ (respectively $b_{3}$). We
also denote representatives of $\varphi(\alpha)$,
$\varphi(\alpha')$, $\varphi(\alpha_{1})$, $\varphi(\alpha_{2})$,
$\varphi(\alpha_{3})$ by $\varphi(a)$, $\varphi(a')$,
$\varphi(a_{1})$, $\varphi(a_{2})$, $\varphi(a_{3})$,
respectively. If, say, $a_1$ is trivial, then $a_1$ bounds either
a disc with a puncture or a M\"obius band. In the first case, we
think of $a_1$ as the puncture its bounds. In the second case, we
think of $a_1$ as the core of the M\"obius band it bounds. ${\rm
Aut } \, C(N)$ has a well defined action on the isotopy classes of
these trivial simple closed curves.

We now examine each of the thirteen cases.

\begin{figure}[hbt]
\begin{center}
\includegraphics[width=12cm]{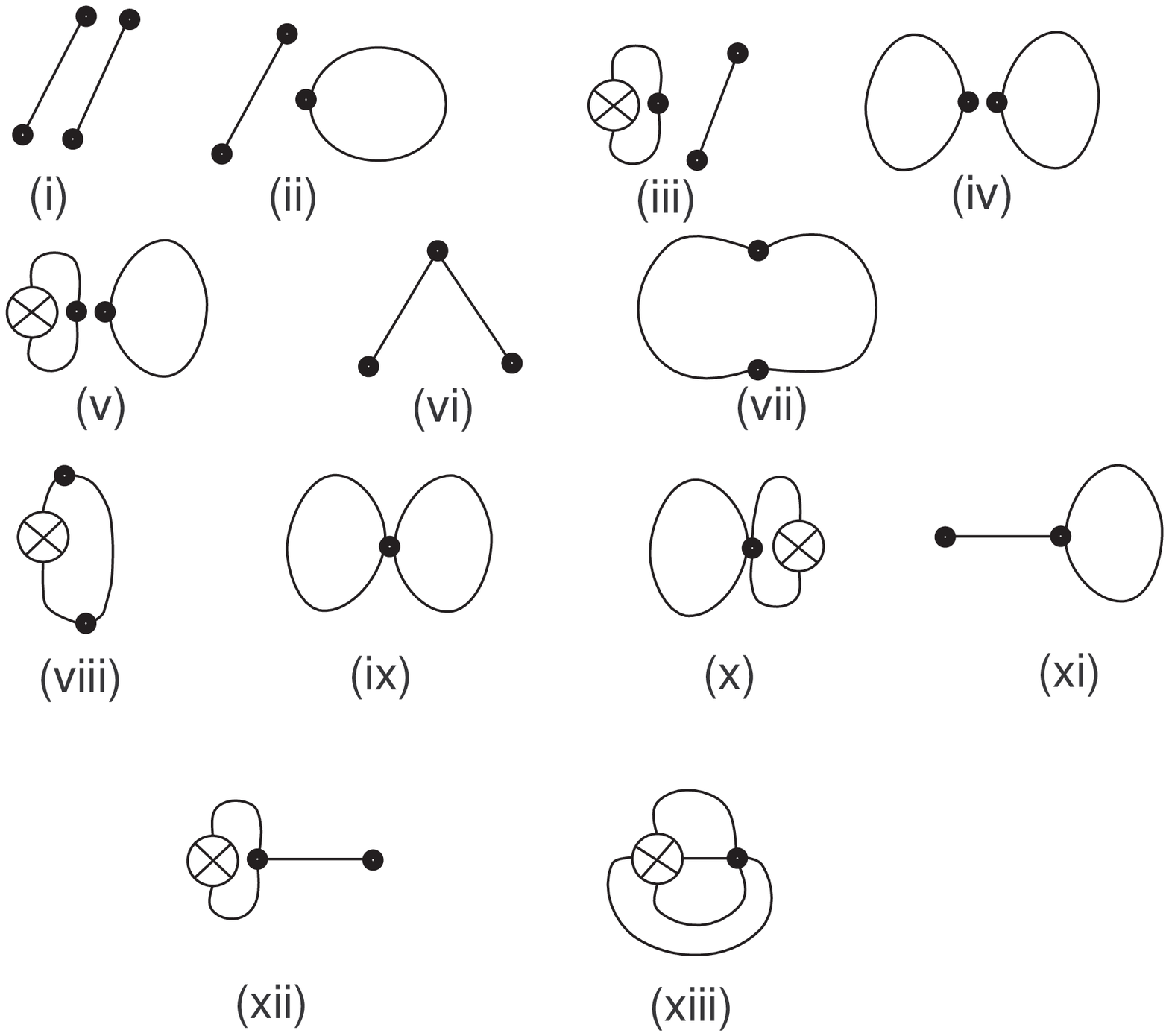}
\caption {} \label{puncture}
\end{center}
\end{figure}

The proof is similar to that of Lemma\,3.7 in \cite{K1} for (i),
(ii), (iv),  (vii), (ix), (xi) and we will not repeat them here.

(iii) In this case, $\varphi(a_{3})$ bounds a M\"obius band $M$
with a puncture and $\varphi(b)$ bounds a disc $D$ with two
punctures. Since $\varphi$ is an automorphism, and since $a_{3}$
and $b$ are disjoint and nonisotopic, $\varphi(a_{3})$ and
$\varphi(b)$ are disjoint and nonisotopic. Then $M$ does not
intersect $D$. Since $\varphi(a')$ is on $M$ and $\varphi(b')$ is
on $D$, it follows that $\varphi(a')$ is disjoint from
$\varphi(b')$.

(v) The once-punctured annulus bounded by $b_{1}$ and $b_{2}$ and
the once-punctured M\"obius band bounded by $a_{3}$ are disjoint.
Since $b_{1}$, $b_{2}$ and $a_{3}$ are pairwise disjoint, so are
$\varphi(b_{1})$, $\varphi(b_{2})$ and $\varphi(a_{3})$. So the
annulus $A$ bounded by $\varphi(b_{1})$, $\varphi(b_{2})$, and the
once-punctured M\"obius band $M$ bounded by $\varphi(a_{3})$ are
disjoint. Since $\varphi(a')$ is on $M$ and $\varphi(b')$ is on
$A$, they are disjoint.

(vi) This case follows from Corollary\,\ref{Corollary-chain}.

(viii) Suppose that $a'$ and $b'$ are joining the punctures $P$
and $Q$. Let $R$ be any other puncture and let $c'$ be an arc from
$P$ to $R$ disjoint from $a' \cup b'$. Let $D$ denote a regular
neighborhood of $b' \cup c' \cup \{P,Q,R\}$, so that $D$ is a disc
with three punctures. Let $d'$ denote the unique arc on $D$
joining $Q$ to $R$ such that $d'$ does not intersect $a' \cup b'
\cup c'$. Let $e$ be the boundary of $D$. So any two arcs in the
set $\{b', c', d' \}$ is a simple pair. Thus any two arcs in the
set $\{\varphi(b'), \varphi(c'), \varphi(d')\}$ is a simple pair,
and $\varphi(b')$, $\varphi(c')$ and $\varphi(d')$ are contained
in the three punctured disk component of $\varphi(e)$. It follows
that any arc disjoint from $\varphi(c')$ and $\varphi(d')$ is also
disjoint from $\varphi(b')$. Since $\varphi(a')$ is disjoint from
$\varphi(c')$ and $\varphi(d')$ by (vi), it is also disjoint from
$\varphi(b')$.

(x) Let $P$ be the common endpoint of two-sided loop $a'$ and
one-sided loop $b'$, so that both arcs connect $P$ to itself. By
Theorem~\ref{Theorem-prese} and Lemma~\ref{Lemma-punc},
$\varphi(a_{1})$ and $\varphi(a_{2})$ are boundaries of an annulus
with one puncture $\varphi(P)$. Since $\varphi(b')$ is disjoint
from $\varphi(a_{2})$, $\varphi(b')$ is also disjoint from
$\varphi(a')$.

(xii) Assume that $a'$ is joining $P$ to itself such that $a'$ is
a one-sided loop, and $b'$ is connecting $P$ to $Q$. Let
$P_{1},...,P_{n-2}$ be the punctures other than $P$ and $Q$.
Choose a chain $\langle c'_{1},...,c'_{n-2} \rangle$ such that
$c'_{i}$ joins $P_{i-1}$ to $P_{i}$ for $1 \leqslant i \leqslant
n-2$, where $P_{0} = Q$. We consider a two-sided loop $d'$ joining
$P$ to $P$ disjoint from $a'$, $b'$ and $\langle
c'_{1},...,c'_{n-2} \rangle$ such that one of the components of
the complement of $d'$ is a M\"obius band and the other is a disc
with punctures $Q, P_{1}, P_{2},...,P_{n-2}$. By (vi), $\langle
\varphi(c'_{1}),...,\varphi(c'_{n-2}) \rangle$ is also chain and
disjoint from $\varphi(b')$. By (xi), $\varphi(d')$ is disjoint
from $\varphi(b')$ and the chain $\langle
\varphi(c'_{1}),...,\varphi(c'_{n-2}) \rangle$. Note that one of
the components of the complement of $\varphi(d')$ contains
$\varphi(b')$ and the chain $\langle
\varphi(c'_{1}),...,\varphi(c'_{n-2}) \rangle$. Since
$\varphi(a')$ is disjoint from $\varphi(d')$ by (x), we obtain
that $\varphi(a')$ is also disjoint from $\varphi(b')$.

(xiii) Let $P$ be the common endpoint one-sided loops $a'$ and
$b'$. The complement of a regular neighborhood of $a'\cup b'$ is
the union of two discs $D_{1}$ and $D_{2}$ with $n-k-1$ and $k$
punctures for some $k$ with $1 \leqslant k \leqslant n-2$. Let
$b'_{1},..., b'_{n-k-2}$ be a chain on $D_{1}$ so that each
$b'_{i}$ is disjoint from $a' \cup b'$. Let $P_{1},...,P_{k}$ be
the punctures on $D_{2}$. We can choose pairwise disjoint arcs
$c'_{i}$ connecting $P_{i-1}$ and $P_{i}$ such that each $c'_{i}$
is also disjoint from $a' \cup b' \cup \partial D_{1}$, where
$P_{0}=P_{k+1}=P$ and $1 \leqslant i \leqslant k+1$.

It follows that a regular neighborhood of $c'_{1} \cup \cdots \cup
c'_{k+1}$ is a M\"obius band with $k+1$ punctures. Then a regular
neighborhood of $\varphi(c'_{1}) \cup \cdots \cup
\varphi(c'_{k+1})$ is also a M\"obius band with $k+1$ punctures.

Now the surface obtained from $N$ by cutting along
$\varphi(c'_{1}) \cup \cdots \cup \varphi(c'_{k+1})$ and
$\varphi(\partial D_{1})$ is an annulus $A$. The puncture
$\varphi(P)$ gives rise to two punctures $R_{1}$ and $R_{2}$ on
the same component of $\partial A$. The arcs $\varphi(a')$ and
$\varphi(b')$ live on $A$. In order to get a one-sided arc, each
must connect $R_{1}$ to $R_{2}$. Up to isotopy there are two arcs
from $R_{1}$ and $R_{2}$ which are disjoint. These two arcs must
be $\varphi(a')$ and $\varphi(b')$.

\end{proof}

\begin{Lemma}\label{Lemma-injective}
The natural map ${\rm Aut} \,C(N) \rightarrow {\rm Aut} \,B(N)$ is
injective.
\end{Lemma}

\begin{proof}
From the lemma above, every element of $\textmd{Aut} \,C(N)$
yields an element of $\textmd{Aut} \,B(N)$. This gives a
homomorphism. Now, we need to show that the kernel of this
homomorphism is trivial. In other words, if an automorphism of
$C(N)$ induces the identity automorphism of $B(N)$, then this
automorphism must be the identity.

Let $\varphi$ be an automorphism of $C(N)$. Suppose that $\varphi$
induces the identity automorphism of $B(N)$. We recall that there
is a one-to-one correspondence between $2$-separating vertices of
$C(N)$ and the vertices of $B(N)$ joining different punctures. It
follows that $\varphi$ is the identity on $2$-separating vertices
of $C(N)$.

Let $\alpha$ be a one-sided vertex of $C(N)$ and let $a \in
\alpha$. Let us denote by $P_{1},\ldots,P_{n}$ the punctures on
the connected component of $N_{a}$. Let us take any chain $\langle
c_{1}',\ldots,c_{n-1}' \rangle$ disjoint from $a$ such that
$c_{i}'$ connects $P_{i}$ to $P_{i+1}$. Let $\gamma_{i}'$ be the
isotopy class of $c_{i}'$. Since $i(\gamma_{i},\alpha) = 0$, we
have $i(\gamma_{i},\varphi(\alpha)) = 0$ and hence
$i(\gamma_{i}',\varphi(\alpha)) = 0$. Let $C' = c_{1}' \cup \cdots
\cup c_{n-1}'$. The surface $N_{C'}$ obtained from $N$ by cutting
along $C'$ is a projective plane with one boundary component. Up
to isotopy, there is only one one-sided simple closed curve on
$N_{C'}$. Both $\alpha$ and $\varphi(\alpha)$ are on $N_{C'}$. So,
we must have $\varphi(\alpha) = \alpha$.

Let $\alpha$ be a $k$-separating vertex of $C(N)$ with $3
\leqslant k \leqslant n-1$ and let $a \in \alpha$. Let us denote
by $P_{1},\ldots,P_{k}$ and $Q_{1},\ldots,Q_{n-k}$ the punctures
on the two connected components of $N_{a}$. Let us take any two
chains $\langle b_{1}',\ldots,b_{k-1}' \rangle$ and $\langle
c_{1}',\ldots,c_{n-k-1}' \rangle$ disjoint from $a$ such that
$b_{i}'$ connects $P_{i}$ to $P_{i+1}$ and $c_{j}'$ connects
$Q_{j}$ to $Q_{j+1}$. Let $\beta_{i}'$ and $\gamma_{j}'$ be the
isotopy classes of $b_{i}'$ and $c_{j}'$, respectively. Let $d'$
be a one-sided loop joining to $Q_{1}$ to itself disjoint from the
other arcs. Let $\delta'$ be the isotopy class of $d'$. Because
$i(\beta_{i},\alpha) = 0$, $i(\gamma_{j},\alpha) = 0$ and
$i(\delta,\alpha) = 0$, we have $i(\beta_{i},\varphi(\alpha)) =
i(\gamma_{j},\varphi(\alpha)) = 0$ and $i(\delta,\varphi(\alpha))
= 0$. Then $i(\beta_{i}',\varphi(\alpha)) =
i(\gamma_{j}',\varphi(\alpha)) = 0$ and
$i(\delta',\varphi(\alpha)) = 0$. Let $A' = b_{1}' \cup \cdots
\cup b_{k-1}' \cup d' \cup c_{1}' \cup \cdots \cup c_{n-k-1}'$.
The surface $N_{A'}$ obtained from $N$ by cutting along $A'$ is an
annulus. Since, up to isotopy, there is only one two-sided simple
closed curve on $N_{A'}$, we must have $\varphi(\alpha) = \alpha$.

\end{proof}

\bigskip
\subsubsection{Ideal triangulations of $N$ and maximal
simplices of $B(N)$}

All maximal simplices in $B(N)$ have the same dimension, and there
is a well-defined action of the group $\textmd{Aut} \,B(N)$ on
maximal simplices. Any realization of a maximal simplex is an
ideal triangulation of $N$. An ideal triangulation is a
triangulation of $N$ whose vertex set is the set of punctures on
$N$ in the sense that vertices of a triangle can coincide as can a
pair of edges. Note that isotopy class of any ideal triangulation
forms a maximal simplex in $B(N)$. The converse of that is also
true. So, $\textmd{Aut} \,B(N)$ acts on the isotopy classes of
ideal triangulations.

We  quote the following definition from \cite{K1}.

\begin{Definition}
A good ideal triangle is a set $\{a',b',c'\}$ of nontrivial
embedded disjoint arcs such that $a'$, $b'$ and $c'$ connect
$P_{1}$ to $P_{2}$, $P_{2}$ to $P_{3}$ and $P_{3}$ to $P_{1}$,
respectively, for three different punctures $P_{1}$, $P_{2}$ and
$P_{3}$, and such that $a' \cup b' \cup c'$ bounds a disc in $N$.
\end{Definition}

The following lemma can be proven similar to the proof Corollary
in \cite{Ha}.

\begin{Lemma}\label{Lemma-Hatcher}
Let $N$ be a projective plane with at least five punctures. Then
given any two maximal simplices $\sigma$ and $\sigma'$ of $B(N)$,
there exists a sequence of maximal simplices
$\sigma=\sigma_{1},\sigma_{2},\ldots,\sigma_{k}=\sigma'$ such that
$\sigma_{i} \cap \sigma_{i+1}$ is a codimension-one simplex for
each $i$.
\end{Lemma}

The lemmas below are analogous to Lemma\,3.9 and Lemma\,3.10 in
\cite{K1}.

\begin{Lemma}\label{Lemma-ideal}
Let $\tilde{\varphi} \in {\rm Aut} \,B(N)$,
$\triangle=\{a',b',c'\}$ be a good ideal triangle and let
$\alpha'$, $\beta'$, $\gamma'$ be the isotopy classes of $a'$,
$b'$, $c'$, respectively. Then $\{\alpha',\beta',\gamma'\}$ and,
hence
$\{\tilde{\varphi}(\alpha'),\tilde{\varphi}(\beta'),\tilde{\varphi}(\gamma')\}$
is a $2$-simplex in $B(N)$. If
$\tilde{\varphi}(\triangle)=\{\tilde{\varphi}(a'),\tilde{\varphi}(b'),\tilde{\varphi}(c')\}$
is a realization of the latter simplex, then it is a good ideal
triangle on $N$.
\end{Lemma}

\begin{Lemma}\label{Lemma-agree}
Let $\tilde{\varphi}$ and $\tilde{\psi}$ be two automorphisms of
$B(N)$. If they agree on a maximal simplex, then they agree on all
of $B(N)$.
\end{Lemma}

\bigskip
\subsubsection{Proof of Theorem\,\ref{Theorem-main1} for punctured
$\mathbb{R}P^{2}$} In Section\,3, we showed that the natural map
$\mathcal M_N \rightarrow \textmd{Aut} \,C(N)$ is injective. We
show that this natural homomorphism is onto. Let $\varphi \in
\textmd{Aut} \,C(N)$ and let $\sigma$ be the isotopy class of a
good ideal triangulation of $N$. So, $\sigma$ is a maximal simplex
of $B(N)$. By Lemma\,\ref{Lemma-ideal}, $\tilde{\varphi} \in
\textmd{Aut} \,B(N)$, the automorphism induced by $\varphi$, takes
a good ideal triangle to a good ideal triangle and
$\tilde{\varphi}$ can be realized by a homeomorphism. Also,
because each edge of a good ideal triangulation is an edge of two
good ideal triangles, the homeomorphism of these triangles gives a
homeomorphism $\Phi$ of $N$. By replacing $\Phi$ by a
diffeomorphism isotopic to $\Phi$ if necessary, we may assume that
$\Phi : N \rightarrow N$ is a diffeomorphism. If $[\Phi]$ is the
isotopy class of $\Phi$, then $\tilde{\varphi}$ agrees with
$\tilde{\Phi}_{*}$ the automorphism induced by $\Phi$, on the
maximal simplex $\sigma$ of $B(N)$. From Lemma\,\ref{Lemma-agree},
they agree on all of $B(N)$. Thus, $\tilde{\varphi} =
\tilde{\Phi}_{*}$. Since the map $\textmd{Aut} \,C(N) \rightarrow
\textmd{Aut} \,B(N)$ is injective, we get $\varphi = \Phi_{*}$.

The proof of the theorem for punctured $\mathbb{R}P^{2}$ is now
complete.

\section{Surfaces of Higher Genus}
\bigskip

Throughout this section unless otherwise stated, $N$ will denote a
connected nonorientable surface of genus $g$ with $n$ punctures
where $g$ is odd and $g+n \geqslant 6$.

In this section, we first show that automorphisms of $C(N)$
preserve the topological type of the vertices of $C(N)$. We then
prove that every automorphism of $C(N)$ is induced by a
diffeomorphism of $N$. We prove this by induction on $r$, where
$g=2r+1$.

\begin{Lemma}\label{Lemma-intersect}
Let $\varphi$ be an automorphism of $C(N)$. Let $\alpha$ and
$\beta$ be nonseparating two-sided vertices. If
$i(\alpha,\beta)=1$, then $i(\varphi(\alpha),\varphi(\beta))=1$.
\end{Lemma}

\begin{proof}
Let $\gamma_{1}$, $\gamma_{2}$,..., $\gamma_{2r-1}$ be pairwise
disjoint one-sided vertices such that each $\gamma_{i}$ is
disjoint from $\alpha$ and $\beta$. Consider the link
$L(\gamma_{1}, ... ,\gamma_{2r-1})$ of these vertices. Then
$\varphi$ restricts to an automorphism $\varphi_{|} :
L(\gamma_{1}, ... ,\gamma_{2r-1}) \rightarrow L(\gamma_{1}, ...
,\gamma_{2r-1})$. $L(\gamma_{1}, ... ,\gamma_{2r-1})$ is
isomorphic to the complex of curves of a torus with $n+2r-1$
punctures. Observe that $\alpha$ and $\beta$ vertices are
contained in $L(\gamma_{1}, ... ,\gamma_{2r-1})$ and by assumption
$i(\alpha,\beta)=1$. Since $n+2r-1 \geqslant 3$, by Theorem\,1 in
\cite{K1}. $\varphi_{|} $ is induced by a diffeomorphism. In
particular,  $i(\varphi(\alpha),\varphi(\beta))=1$.
\end{proof}

The following lemma can be proven using techniques similar to
those in Lemma\,3.8. in \cite{Ir2} or Lemma\,3 in \cite{IK}.

\begin{Lemma}\label{Lemma-sequence}
Let $a$ be a two-sided nonseparating circle. If $c$ and $d$ are
two two-sided nonseparating circles, both intersecting $a$
transversely once, then there is a sequence $c=c_{0},
c_{1},...,c_{n}=d$ of two-sided nonseparating circles such that
each $c_{i}$ intersects the circle $a$ transversely once and
$c_{i}$ is disjoint from $c_{i+1}$ for each $i=0,1,...,n-1$.
\end{Lemma}

\begin{Lemma}\label{Lemma-notiso}
Let $n \geqslant 2$ and $k \geqslant 1$. If $N$ is a projective
plane with $n$ punctures and $T$ is a torus with $k$ punctures,
then $C(N)$ and $C(T)$ are not isomorphic.
\end{Lemma}

\begin{proof}
The complexes $C(N)$ and $C(T)$ have dimensions $n-2$ and $k-1$,
respectively. If $k \neq n-1$, since these complexes of curves
have different dimensions, $C(N)$ and $C(T)$ are not isomorphic.
If $k = n-1$, we proceed as follows.

Let $n=2$. Thus $N$ is a projective plane with two punctures and
$T$ is a torus with one puncture. In this case, $C(N)$ is finite
(\cite{S}), however, $C(T)$ is infinite discrete since $T$ is a
torus with one puncture. Therefore, they are not isomorphic.

Now, assume that $C(N)$ and $C(T)$ are not isomorphic when $N$ has
$n-1$ punctures. Since there are $n-1$ punctures on $N$, there are
$n-2$ punctures on $T$. We need to show that the complexes $C(N)$
and $C(T)$ are not isomorphic if $N$ has $n$ punctures. Assume the
contrary that there is an isomorphism $\varphi : C(N) \rightarrow
C(T)$. For any vertex $\gamma$ in $C(N)$, $\varphi$ induces an
isomorphisms $L(\gamma) \rightarrow L(\varphi(\gamma))$ and
$L^{d}(\gamma) \rightarrow L^{d}(\varphi(\gamma))$.

Note that for a vertex $\gamma$ of $C(N)$, the dual link
$L^{d}(\gamma)$ of $\gamma$ is connected if and only if $\gamma$
is either one-sided or $2$-separating. For a vertex $\delta$ of
$C(T)$, the dual link $L^{d}(\delta)$ of $\delta$ is connected if
and only if $\delta$ is either nonseparating or $2$-separating
vertex. From this, it follows that the image of the union of the
set of one-sided vertices and the set of $2$-separating vertices
of $C(N)$ is precisely the union of the set of nonseparating
vertices and the set of $2$-separating vertices of $C(T)$. Let
$\gamma$ be a $2$-separating vertex of $C(N)$. Then
$\varphi(\gamma)$ is either a nonseparating vertex or a
$2$-separating vertex of $C(T)$. First, we assume that
$\varphi(\gamma)$ is a nonseparating vertex $\delta$ of $C(T)$.
Clearly, $L(\gamma)$ is isomorphic to the complex of curves of a
projective plane with $n-1$ punctures and $L(\delta)$ is
isomorphic to the complex of curves of a sphere with $n+1$
punctures. By Lemma~\ref{Lemma-iso}, these complexes are not
isomorphic. Therefore, $\varphi(\gamma)$ cannot be a nonseparating
vertex $\delta$ of $C(T)$. Now, we assume that $\varphi(\gamma)$
is a $2$-separating vertex $\mu$ of $C(T)$. Then $L(\mu)$ is
isomorphic to the complex of curves of a torus with $n-2$
punctures. By assumption, these complexes are not isomorphic.
Therefore, we get a contradiction. Hence, $C(N)$ and $C(T)$ are
not isomorphic.

\end{proof}

\begin{Lemma}\label{Lemma-notisogenel}
Let $r \geqslant 1$ and $n \geqslant 0$. If $N$ is a connected
nonorientable surface of genus $2r+1$ with $n$ punctures and $M$
is a connected nonorientable surface of genus $2r$ with $k$
punctures, then $C(N)$ and $C(M)$ are not isomorphic.
\end{Lemma}

\begin{proof}
The complexes $C(N)$ and $C(M)$ have dimensions $4r+n-2$ and
$4r+k-4$, respectively. If $k \neq n+2$, since these complexes of
curves have different dimensions, $C(N)$ and $C(M)$ are not
isomorphic. If $k = n+2$, we proceed as follows.

Let $s=r+n$. We will prove that $C(N)$ and $C(M)$ are not
isomorphic by induction on $s$. Let $s=1$. In this case, $r=1$ and
$n=0$. In other words, $N$ is a connected closed nonorientable
surface of genus $3$ and $M$ is a Klein bottle with $2$ punctures.
We need to show that the complexes $C(N)$ and $C(M)$ are not
isomorphic. Assume the contrary that there is an isomorphism
$\varphi : C(N) \rightarrow C(M)$. For any vertex $\gamma$ in
$C(N)$, $\varphi$ induces an isomorphisms $L(\gamma) \rightarrow
L(\varphi(\gamma))$ and $L^{d}(\gamma) \rightarrow
L^{d}(\varphi(\gamma))$. Note that for a vertex $\gamma$ of
$C(N)$, the dual link $L^{d}(\gamma)$ of $\gamma$ is connected if
and only if $\gamma$ is either one-sided vertex and let $c \in
\gamma$ such that $N_{c}$ is an orientable surface or one-sided
vertex such that $N_{c}$ is a nonorientable surface or
nonseparating two-sided vertex. For a vertex $\delta$ of $C(M)$,
the dual link $L^{d}(\delta)$ of $\delta$ is connected if and only
if $\delta$ is either one-sided or nonseparating two-sided or
$2$-separating vertex. Let $\gamma$ be a nonseparating two-sided
vertex of $C(N)$. Then $\varphi(\gamma)$ is either a one-sided
vertex or a nonseparating two-sided vertex or a $2$-separating
vertex of $C(M)$. First, we assume that $\varphi(\gamma)$ is a
one-sided vertex $\delta$ of $C(M)$. Obviously, $L(\gamma)$ is
isomorphic to the complex of curves of a projective plane with $2$
punctures and $L(\delta)$ is isomorphic to the complex of curves
of a projective plane with $3$ punctures. Since these complexes
have different dimensions zero and one, respectively, these
complexes are not isomorphic. Second, we suppose that
$\varphi(\gamma)$ is a nonseparating two-sided vertex $\mu$ of
$C(M)$. Then $L(\mu)$ is isomorphic to the complex of curves of a
sphere with $4$ punctures. By Lemma~\ref{Lemma-iso}, these
complexes are not isomorphic. Now, we assume that
$\varphi(\gamma)$ is a $2$-separating vertex $\lambda$ of $C(M)$.
Clearly, $L(\lambda)$ is isomorphic to the complex of curves of a
Klein bottle with one puncture. Since these complexes have
different dimensions, these complexes are not isomorphic.
Therefore, we get a contradiction. Hence, the complex of curves of
a connected closed nonorientable surface of genus $3$ and the
complex of curves of a Klein bottle with $2$ punctures are not
isomorphic.

Assume that the complexes $C(N)$ and $C(M)$ are not isomorphic for
all integers $1$ to $s-1$. Now, let us prove it for $s$. In this
case, $N$ is a connected nonorientable surface of genus $2r+1$
with $n$ punctures and $M$ is a connected nonorientable surface of
genus $2r$ with $n+2$ punctures. Suppose the contrary that there
is an isomorphism $\varphi : C(N) \rightarrow C(M)$. Note that for
any vertex $\alpha$ of $C(N)$, the dual link $L^{d}(\alpha)$ of
$\alpha$ is connected if and only if $\alpha$ is either one-sided
vertex and let $a \in \alpha$ such that $N_{a}$ is an orientable
surface or one-sided vertex such that $N_{a}$ is a nonorientable
surface or nonseparating two-sided vertex or $2$-separating
vertex. For a vertex $\beta$ of $C(M)$, the dual link
$L^{d}(\beta)$ of $\beta$ is connected if and only if $\beta$ is
either one-sided vertex or nonseparating two-sided vertex and $b
\in \beta$ such that $M_{b}$ is an orientable surface or
nonseparating two-sided vertex such that $M_{b}$ is a
nonorientable surface or $2$-separating vertex. Let $\alpha$ be a
$2$-separating vertex of $C(N)$. First,we assume that
$\varphi(\alpha)$ is a one-sided vertex $\zeta$ of $C(M)$.
Clearly, $L(\alpha)$ is isomorphic to the complex of curves of a
connected nonorientable surface of genus $2r+1$ with $n-1$
punctures and $L(\zeta)$ is isomorphic to the complex of curves of
a connected nonorientable surface of genus $2r-1$ with $n+3$
punctures. Although these complexes have the same dimensions,
these complexes are not isomorphic. Because, by
Proposition~\ref{Prop-ndim1}, there is a maximal simplex of
dimension $3r+n-3$ in the complex of curves $C(N_{2r+1,n-1})$,
whereas there is no any maximal simplex of dimension $3r+n-3$ in
the complex of curves $C(M_{2r-1,n+3})$. Second, we suppose that
$\varphi(\alpha)$ is a nonseparating two-sided vertex $\epsilon$
of $C(M)$ and let $e \in \epsilon$ such that $M_{e}$ is an
orientable surface. Then $L(\epsilon)$ is isomorphic to the
complex of curves of a connected orientable surface of genus $r-1$
with $n+4$ punctures. By Lemma~\ref{Prop-odim1} and
Proposition~\ref{Prop-ndim1}, these complexes $C(M_{r-1,n+4})$ and
$C(N_{2r+1,n-1})$ are not isomorphic. Third, we assume that
$\varphi(\alpha)$ is a nonseparating two-sided vertex $\omega$ of
$C(M)$ and let $w \in \omega$ such that $M_{w}$ is a nonorientable
surface. Then $L(\omega)$ is isomorphic to the complex of curves
of a connected nonorientable surface of genus $2r-2$ with $n+4$
punctures. Since these complexes $C(M_{2r-2,n+4})$ and
$C(N_{2r+1,n-1})$ have different dimensions $4r+n-4$ and $4r+n-3$,
respectively, these complexes are not isomorphic. Now, we suppose
that $\varphi(\alpha)$ is a $2$-separating vertex $\nu$ of $C(M)$.
Then $L(\nu)$ is isomorphic to the complex of curves of a
connected nonorientable surface of genus $2r$ with $n+1$
punctures. By assumption of the induction, these complexes
$C(M_{2r,n+1})$ and $C(N_{2r+1,n-1})$ are not isomorphic.
Therefore, we get a contradiction. Hence, $C(N)$ and $C(M)$ are
not isomorphic.

\end{proof}

\begin{Lemma}\label{Lemma-take}
Let $N$ be a connected nonorientable surface of genus $g=2r+1,$ $r
\geqslant 1$ with $n \geqslant 0$ boundary components. The group
${\rm Aut} \,C(N)$ preserves the topological type of the vertices
of $C(N)$.
\end{Lemma}

\begin{proof}
Let $\varphi$ be an automorphism of $C(N)$. Note that for a vertex
$\gamma$ of $C(N)$, the dual link $L^{d}(\gamma)$ of $\gamma$ is
connected if and only if $\gamma$ is either one-sided vertex or
nonseparating two-sided vertex or $2$-separating. Therefore,
$\varphi$ cannot take a one-sided vertex or a nonseparating
two-sided vertex or a $2$-separating to a $k$-separating vertex
with $k > 2$ or separating vertex.

Let $\alpha$ be a one-sided vertex of $C(N)$ and let $a \in
\alpha$ such that $N_{a}$ is an orientable surface. Clearly,
$L(\alpha)$ is isomorphic to the complex of curves of an
orientable surface of genus $r$ with $n+1$ boundary components
$S_{r,n+1}$. Let $\beta$ be a one-sided vertex and let $b \in
\beta$ such that $N_{b}$ is a nonorientable surface of genus $2r$
with $n+1$ boundary components. Obviously, $L(\beta)$ is
isomorphic to the complex of curves of a nonorientable surface of
genus $2r$ with $n+1$ boundary components $N_{2r,n+1}$. Let
$\gamma$ be a nonseparating two-sided vertex of $C(N)$ and
$L(\gamma)$ is isomorphic to the complex of curves of a
nonorientable surface of genus $2r-1$ with $n+2$ boundary
components $N_{2r-1,n+2}$. Let $\delta$ be a $2$-separating vertex
of $C(N)$ and $L(\delta)$ is isomorphic to the complex of curves
of a nonorientable surface of genus $2r+1$ with $n-1$ boundary
components $N_{2r+1,n-1}$. Since there are maximal simplices of
different dimensions in the complexes of curves $C(N_{2r,n+1})$,
$C(N_{2r-1,n+2})$ and $C(N_{2r+1,n-1})$, these complexes are not
isomorphic to $C(S_{r,n+1})$. Moreover, the dimensions of
$C(N_{2r,n+1})$, $C(N_{2r-1,n+2})$ and $C(N_{2r+1,n-1})$ are
$4r+n-3$, $4r+n-4$ and $4r+n-3$, respectively. Therefore, the
complex of curves $C(N_{2r-1,n+2})$ is not isomorphic to
$C(N_{2r,n+1})$ and $C(N_{2r+1,n-1})$. Furthermore, although these
complexes of curves $C(N_{2r,n+1})$ and $C(N_{2r+1,n-1})$ have the
same dimensions, by Lemma~\ref{Lemma-notisogenel}, these complexes
are not isomorphic.

Let $\lambda$ be a separating vertex and let $e \in \lambda$ be
circle. $N_{e}$ has two connected components. Let $N_{e} = N_{0}
\sqcup N_{1}$ such that $N_{0}$ and $N_{1}$ are nonorientable
surfaces. More precisely, $N_{0}$ is a nonorientable surface of
genus $l$ for some $1\leqslant l \leqslant 2r$ with $k+1$
punctures for some $0 \leqslant k \leqslant n$ denoted by
$N_{l,k+1}$, then $N_{1}$ is a nonorientable surface of genus
$2r+1-l$ with $n-k+1$ punctures denoted by $N_{2r+1-l,n-k+1}$. The
dual link $L^{d}(\lambda)$ has exactly two connected components.
Let us denote these components by $L_{0}^{d}(\lambda)$ and
$L_{1}^{d}(\lambda)$; $L^{d}(\lambda) = L_{0}^{d}(\lambda) \sqcup
L_{1}^{d}(\lambda)$. We name $L_{i}(\lambda)$ so that the vertices
of $L^{d}_{i}(\lambda)$ are the isotopy classes of circles on
$N_{i}$. Let $L_{i}(\lambda)$ be the full subcomplex of $C(N)$
with vertices $L^{d}_{i}(\lambda)$. It follows that
$(L_{i}(\lambda))^{d} = L^{d}_{i}(\lambda)$. Then $L_{i}(\lambda)$
is isomorphic to $C(N_{i})$. Let $\mu$ be a separating vertex and
$w \in \mu$ be circle. Similarly, we define $L^{d}_{i}(\mu)$ and
$L_{i}(\mu)$ so that $L_{i}(\mu)$ is isomorphic to $C(N'_{i})$. In
other words, $L_{0}(\mu)$ is isomorphic to the complex of curves
of an orientable surface of genus $l$ for some $0\leqslant l
\leqslant r$ with $k+1$ punctures for some $0 \leqslant k
\leqslant n$ denoted by $S_{l,k+1}$, then $L_{1}(\mu)$ is
isomorphic to the complex of curves of a nonorientable surface of
genus $2(r-l)+1$ with $n-k+1$ denoted by $N_{2(r-l)+1,n-k+1}$.
Since $L_{0}(\mu)$ is isomorphic to $C(S_{l,k+1})$ and all maximal
simplices in $C(S_{l,k+1})$ have the same dimension, $L^{d}(\mu)$
is not isomorphic to $L^{d}(\lambda)$. Hence, $\lambda$ cannot be
mapped $\mu$ under $\varphi$.

The proof of the lemma is complete.
\end{proof}

\begin{Remark}\label{Remark-take}
In the above proof, in case of $r=1$ and $n \geqslant 0$, one can
see that the complex of curves of a torus with $n+1$ boundary
components and the complex of curves of a projective plane with
$n+2$ boundary components are not isomorphic by
Lemma~\ref{Lemma-notiso}.
\end{Remark}

\begin{Theorem}\label{Thm-onto}
Let $N$ be a connected nonorientable surface of odd genus $g=2r+1$
with $n$ punctures. Suppose that $g+n \geqslant 6$. Then $\varphi$
agrees with a map $h_{*} : C(N) \rightarrow C(N)$ which is induced
by a diffeomorphism $h : N \rightarrow N$.
\end{Theorem}

\begin{proof}

Let $r=0$. Then $N$ is a projective plane with $n \geqslant 5$
punctures. In Section 4, we showed that $\varphi$ is induced by a
diffeomorphism of $N$.

We assume that every automorphism $C(N) \rightarrow C(N)$ is
induced by a diffeomorphism of $N$ if $N$ is of odd genus $g
\leqslant 2r-1$. Now, we show that $\varphi$ is induced by a
diffeomorphism of $N$ for genus $g=2r+1$.

Let $c$ be any two-sided nonseparating circle and $\gamma$ denote
its isotopy class. By Lemma~\ref{Lemma-take}, $\varphi$ takes
$\gamma$ to a two-sided nonseparating vertex, say $\gamma'$, and
let $c' \in \gamma'$. There is a diffeomorphism $f$ such that
$f(c)=c'$. Then $f_{*}^{-1}\varphi(\gamma)=\gamma$ where $f_{*} :
C(N) \rightarrow C(N)$ the automorphism induced by $f$. By
replacing $\varphi$ by $f_{*}^{-1}\varphi$, we can assume that
$\varphi(\gamma)=\gamma$. $\varphi$ restricts to an automorphism
$\varphi_{c} : L(\gamma) \rightarrow L(\gamma)$. Since $L(\gamma)$
is isomorphic to the complex of curves $C(N_{c})$ of a
nonorientable surface $N_{c}$ of genus $2r-1$ with $n+2$ boundary
components, we get an automorphism $\varphi_{c}$ of the complex of
curves of a nonorientable surface of genus $2r-1$ with $n+2$
boundary components. By induction, we can assume that
$\varphi_{c}$ is equal to a map $(\overline{g}_{c})_{*} : C(N_{c})
\rightarrow C(N_{c})$ which is induced by a diffeomorphism
$\overline{g}_{c} : N_{c} \rightarrow N_{c}$. By gluing two
boundary components of $N_{c}$ obtained $c$ in a convenient way
$g_{c}$ induces a diffeomorphism $g_{c} : N \rightarrow N$. It
follows that $g_{c}(c) = c$. Therefore, $\varphi$ agrees with
$(g_{c})_{*}$ on every element of $L(\gamma)$. The composition
$(g^{-1}_{c})_{*}$ with $\varphi$, $(g^{-1}_{c})_{*} \circ
\varphi$ fixes $\gamma$ and every element of $L(\gamma)$. We may
replace $(g^{-1}_{c})_{*} \circ \varphi$ by $\varphi$. Now,
$\varphi$ is an automorphism of $C(N)$ such that it is identity on
$\gamma \cup L(\gamma)$.

Let $d$ be a two-sided nonseparating circle dual to $c$. In other
words, $d$ intersects $c$ transversely once and there is no other
intersection. Let $T$ be a regular neighborhood of $c \cup d$. The
surface $T$ is a torus with one boundary component. Let $e$ be the
boundary component of $T$. Let $\delta$ and $\epsilon$ denote the
isotopy class of $d$ and $e$, respectively. Obviously, since
$\epsilon \in L(\gamma)$, $\varphi(\epsilon) = \epsilon$. Also,
$\varphi(\gamma)=\gamma$. Since $i(\delta,\epsilon)=0$,
$i(\varphi(\delta),\varphi(\epsilon))=i(\varphi(\delta),\epsilon)=0$.
Since $i(\gamma,\delta)=1$ by Lemma~\ref{Lemma-intersect}, there
exists an integer $n$ such that
$\varphi(\delta)=(t_{c}^{n})_{*}(\delta)$.

Let $D(\gamma)$ be the set of isotopy classes of two-sided
nonseparating circles which are dual to $c$ on $N$. Let $d_{1}$ be
a two-sided nonseparating circle which is disjoint from $d$ and
dual to $c$. Similarly, there exists an integer $n_{1}$ such that
$\varphi(\delta_{1})=(t_{c}^{n_{1}})_{*}(\delta_{1})$, where
$\delta_{1}$ is the isotopy class of $d_{1}$. Since
$i(\delta,\delta_{1})=0$,
$i(\varphi(\delta),\varphi(\delta_{1}))=0$. If $n \neq n_{1}$,
then $i((t_{c}^{n})_{*}(\delta), (t_{c}^{n_{1}})_{*}(\delta_{1}))
\neq 0$ since both $d$ and $d_{1}$ are dual to $c$. Thus, we have
$i(\varphi(\delta), \varphi(\delta_{1})) \neq 0$. This is a
contradiction. Therefore, we must have $n=n_{1}$. Now for any
two-sided nonseparating circle $s$ which is dual to $c$, by
Lemma~\ref{Lemma-sequence} we can find a sequence of two-sided
nonseparating circles dual to $c$, connecting $d$ to $s$ such that
each consecutive pair is disjoint. It follows that $\varphi$
agrees with $(t_{c}^{n})_{*}$ on every element of $D(\gamma)$.
Therefore, $\varphi$ agrees with $(t_{c}^{n})_{*}$ on $\gamma$ and
on $L(\gamma)\cup D(\gamma)$. Let $h_{c} = t_{c}^{n}$. In the
following, we will denote $(h_{c})_{*}$ simply by $h_{*}$.

If $u$ and $v$ are any other two-sided nonseparating circles dual
to each other, then there exists a diffeomorphism $N \rightarrow
N$ mapping $(u,v)$ to $(c_{1},c_{2})$, where $c_{1}$ and $c_{2}$
are the circles in Figure\,\ref{sur-circles}. The set
$\{c_{1},c_{2}\}$ can be completed to a set $C$ of two-sided
circles except $a$ as shown in Figure\,\ref{sur-circles}. Let
$\gamma_{i}$ be the isotopy class of $c_{i}$.

We orient tubular neighborhoods of elements of $C$ in such way
that these orientation agree with an orientation of a regular
neighborhood $M$ of $\cup c_{i}$. We note that $M$ is an
orientable surface.

The isotopy classes of Dehn twists about elements of $C$ generate
a subgroup $\mathcal{T'}$. As shown above, the restriction of
$\varphi$ to $L(\gamma_{i}) \cup D(\gamma_{i})$ agrees with the
induced map $(h_{i})_{*}$ of a diffeomorphism $h_{i} : N
\rightarrow N$.

Any circle in $C$ is either disjoint from $c_{1}$ or is dual to
$c_{1}$, $\varphi(\gamma_{i})=(h_{1})_{*}(\gamma_{i})$ for all
$i$. Similarly, $\varphi(\gamma_{i})=(h_{2})_{*}(\gamma_{i})$ for
all $i$. Hence, $(h^{-1}_{1} \circ h_{2})_{*}(\gamma_{i}) =
\gamma_{i}$. In particular, $(h^{-1}_{1} \circ h_{2})_{*}(c_{i})$
is isotopic to $c_{i}$ and
 $$(h_{1}^{-1} \circ h_{2})t_{c_{i}}(h_{1}^{-1} \circ
h_{2})^{-1}=t_{c_{i}}^{\epsilon_{i}},$$ where $\epsilon_{i} =
\pm1$. Let us denote $h_{1}^{-1} \circ h_{2}$ by $h$. For a circle
$c_{i} \in C$ dual to $c_{1}$, we have the braid relation
$$t_{c_{1}}t_{c_{i}}t_{c_{1}}=t_{c_{i}}t_{c_{1}}t_{c_{i}}.$$
By conjugating with $h$, we get
$$t^{\epsilon_{1}}_{c_{1}}t^{\epsilon_{i}}_{c_{i}}t^{\epsilon_{1}}_{c_{1}}
=
t^{\epsilon_{i}}_{c_{i}}t^{\epsilon_{1}}_{c_{1}}t^{\epsilon_{i}}_{c_{i}}.$$
But this relation holds if and only if $\epsilon_{1}$ and
$\epsilon_{i}$ has the same sign. Now if $c_{i}$ is dual to
$c_{1}$ and if $c_{j}$ is dual to $c_{i}$, by similar reasoning we
get $\epsilon_{j}$ and $\epsilon_{1}$ have the same sign. It
follows that all $\epsilon_{i}$ have the same sign, as we can pass
from $c_{1}$ to any circle in $C$ through dual circles.

\begin{figure}[hbt]
\begin{center}
\includegraphics[width=12cm]{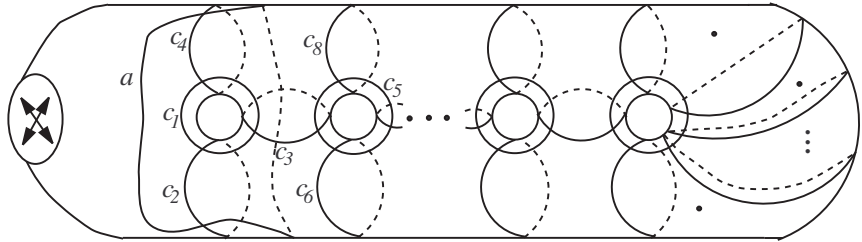}
\caption {Identify antipodal points on the boundary component on
the left hand side} \label{sur-circles}
\end{center}
\end{figure}

Suppose that $\epsilon_{i} = -1$, so that $ht_{c_{i}}h^{-1} =
t^{-1}_{c_{i}}$ for all $c_{i} \in C$. Thus $h$ reverses the
orientation of tubular neighborhoods of $c_{i}$. Let $\rho$ be the
reflection in $yz-$plane as shown in Figure\,\ref{sur-circles}.
The reflection $\rho$ leaves each $c_{i} \in C$ invariant, but
reverses the orientations of tubular neighborhoods of $c_{i}$.
Here $\rho t_{c_{i}} \rho^{-1} = t^{-1}_{c_{i}}$ for all $i$. Thus
$$(\rho \circ h)t_{c_{i}}(\rho \circ h)^{-1} = t_{c_{i}}.$$
In other words, $\rho \circ h \in \mathcal{C}_{\mathcal
M_{N}}(\mathcal{T'})$. Since $\mathcal{C}_{\mathcal
M_{N}}(\mathcal{T'})$ is trivial by
Proposition~\ref{Prop-centralizer}, we get that $$h = \rho^{-1} =
\rho.$$ On the other hand, there exists a two-sided nonseparating
circle $a$ as shown in Figure\,\ref{sur-circles} such that the
class $\alpha$ of $a$ is contained in $L(\gamma_{1}) \cap
D(\gamma_{2})$ and $h(\alpha) = \rho(\alpha) \neq \alpha$.
However, this is a contradiction since $h(\alpha) =
\varphi(\alpha) = \alpha$.

Therefore, we must have $\epsilon_{i} = 1$ for all $i$. Thus,
$ht_{c_{i}}h = t_{c_{i}}$ for all $i$. In particular $h =
h_{1}^{-1} \circ h_{2} \in \mathcal{C}_{\mathcal
M_{N}}(\mathcal{T'}) = 1$. Consequently, $h_{1} = h_{2}$.

Let $d$ be any two-sided nonseparating circle and $\delta$ denote
the isotopy class of $d$ on $N$. Since $N$ is a nonorientable
surface of genus $g \geqslant 3$, $c_{1}$ and $d$ are dually
equivalent by Theorem\,3.1 of \cite{K2}. In other words, there
exists a sequence of two-sided nonseparating circles
$a_{1},\ldots,a_{k}$ on $N$ such that $a_{1}=c_{1}$, $a_{k}=d$ and
the circles $a_{i}$ and $a_{i+1}$ are dual. Using this sequence,
we obtain that $(h_{1})_{*}=(h')_{*}=\varphi$. Indeed, for any
two-sided nonseparating circle $d$, $(h_{1})_{*}$ agrees with
$\varphi$. Furthermore, since the isotopy classes of every
separating circle and of any one-sided circle are in $L(\zeta)$
where $\zeta$ is the isotopy class of some two-sided nonseparating
circle $s$, we obtain that $(h_{1})_{*}$ agrees with $\varphi$ on
$C(N)$.
\end{proof}

\bigskip
\subsubsection{Proof of Theorem\,\ref{Theorem-main1} for nonorientable
surface of odd genus} It is shown that $\mathcal M_N \rightarrow
\textmd{Aut} \,C(N)$ is injective in Section\,3.
Theorem\,\ref{Thm-onto} implies that automorphisms of $C(N)$ are
induced by diffeomorphisms of $N$. Hence, this completes the proof
Theorem\,\ref{Theorem-main1}.

\bigskip
\bigskip
\providecommand{\bysame}{\leavevmode\hboxto3em{\hrulefill}\thinspace}

\end{document}